\documentclass[a4paper]{article}
\usepackage{latexsym}
\usepackage{amssymb}
\begin{document}
\title{Proposing $(\omega _1 ,\beta )$-morasses for $\omega _1 \leq \beta$}
\author{Bernhard Irrgang}
\date{December 16, 2010}
\maketitle  
\begin{abstract}
Firstly, I propose a notion of $(\omega _1 , \beta )$-morass for the case that $\omega _1 \leq \beta$. Secondly, I define $\kappa$-standard morasses such that every $\omega _{1+\beta}$-standard morass is an $(\omega _1, \beta)$-morass. Thirdly, I justify these notions by proving: If there is a $\kappa$-standard morass, then there is an $L_\kappa [X]$ with $Card^{L_\kappa [X]}=Card \cap \kappa$ for which the fine structure theory and condensation hold.
\end{abstract}
\section{Introduction}
In set theory, structures are often obtained by first recursively constructing small structures and then taking a direct limit to get a bigger one. Usually a chain of structures of size $<$ $\kappa$ is constructed by induction along a cardinal $\kappa$. In this way, a direct limit of size $\kappa$ can be obtained. Morasses are index sets along which structures of size $<$ $\omega_\alpha$ can be constructed by induction in  such a way that the limit has size $\omega_{\alpha+\beta}$. The appropriate morass for such a construction is called an $(\omega_\alpha,\beta)$-morass.
\smallskip\\
Morasses were invented by R. Jensen in the early 1970s. He used them to prove the model-theoretic cardinal transfer theorems (see [ChKe]) in G\"odel's constructible universe $L$. If $\kappa, \lambda$ are infinite cardinals, a structure $\frak{A}$ is said to have type $(\kappa,\lambda)$ if $\frak{A}=\langle A, X^\frak{A},\dots\rangle$ where $card(A)=\kappa$ and $card(X^\frak{A})=\lambda$. The simplest cardinal transfer theorem states that if $\frak{A}$ is a structure of type $(\kappa^+,\kappa)$ then there exists a structure $\frak{B}$ of the same language of type $(\omega_1,\omega)$ which is elementary equivalent to $\frak{A}$. This is proved by the construction of an elementary chain that has $\frak{B}$ as its direct limit. Using morasses, Jensen obtains in $L$ statements of this type for bigger gaps between $\kappa$ and $\lambda$.
\smallskip\\
The most general cardinal transfer theorem is shown in a hand-written set of notes [Jen]. Here, he defines the notion of $(\kappa , \beta )$-morasses for $\beta <\kappa$, $\kappa$ regular. How an $(\omega _1,1)$-morass is used to prove the gap-2 cardinal transfer theorem may be found in [Dev]. The theory of morasses is very far developed and very well examined. In particular it is known how to construct morasses in $L$ (see [Dev], [Fri], [Jen]) and how to force them ([Sta1], [Sta2]). Moreover, D. Velleman has defined so-called simplified morasses, along which morass constructions can be carried out very easily compared to classical morasses  ([Vel1], [Vel2], [Vel4]). They are equivalent to usual morasses ([Don], [Mor]). Besides the cardinal transfer theorems, there are many combinatorial applications of morasses. One is for example the construction of $\kappa ^{++}$-super-Suslin trees by S. Shelah and L. Stanley [ShSt1]. Other applications need strengthenings of morasses, like simplified morasses with linear limits [Vel4].
\smallskip\\
For the case $\kappa \leq \beta$, $(\kappa ,\beta )$-morasses have never been defined. I want to propose a notion of $(\omega _1 , \beta )$-morass for this case. In addition, I will define $\kappa$-standard morasses such that every $\omega _{1+\beta}$-standard morass is an $(\omega _1, \beta)$-morass. I will prove that if there is a $\kappa$-standard morass, then there is an $L_\kappa [X]$ with $Card^{L_\kappa [X]}=Card \cap \kappa$ for which the fine structure theory and condensation hold. In a forthcoming paper [Irr2], I show that if there is an $L_\kappa [X]$ with $Card^{L_\kappa [X]}=Card \cap \kappa$ for which the fine structure theory and condensation hold, then there is a $\kappa$-standard morass. On the one hand, this justifies my definitions. On the other hand, it shows that the definition of $\kappa$-standard morass is best possible in the sense that it completely captures the combinatorics of an  $L_\kappa [X]$ with $Card^{L_\kappa [X]}=Card \cap \kappa$, fine structure theory and condensation. Moreover, I conjecture that the existence of an $\omega _{1+\beta}$-standard morass is actually equivalent to the existence of an $(\omega _1,\beta)$-morass.
\smallskip\\
One notion that is related to my definitions of $(\omega _1, \beta)$-morass for $\omega _1 \leq \beta$ and $\kappa$-standard morass is the premorass that was studied by Jensen in the context of his proofs of global square. In [DJS], H.-D. Donder, R. Jensen and L. Stanley derive from the existence of an appropriate premorass that global square and the combinatorial principle squared scales holds. But they derive global square and squared scales directly from the premorass they construct in $L$ without explicitly axiomatizing the notion of premorass. A similiar approach is followed in [BJW] to provide the necessary combinatorics for the proof of Jensen's coding theorem. Squared scales was formulated by Avraham and Shelah for their work on strong covering (see [She], chapter VIII]). A strengthening of squared scales, which is also proved by the same approach [ShSt2], was used in [ShSt3] by S. Shelah and L. Stanley to give a combinatorial proof of Jensen's coding theorem.
\smallskip\\
Another related notion is that of a smooth category which was introduced by R. Jensen and M. Zeman to prove global square in the core model for measures of order 0 [JeZe]. Similiar systems were studied (again without giving an axio\-matic account) in [SchZe1] and [ScheZe2] by E. Schimmerling and M. Zeman to prove that Jensen extender models satisfy the Gap 1 Morass principle and $\Box _\kappa$ for all $\kappa$ that are not subcompact.
\smallskip\\
It is a natural question to ask in which inner models $(\omega _1, \beta)$-morasses and $\kappa$-standard morasses exist. By the usual argument that $\omega _1$-Erd\"os cardinals do not exist in $L$ (see e.g. theorem V 1.8 of [Dev]), it is easy to see that an inner model $M$ with an $\omega _1$-Erd\"os cardinal cannot be of the form $M=L[X]$ such that $L[X]$ satisfies condensation. But that does not mean, that it is impossible that inner models with $\omega _1$-Erd\"os cardinals (or even larger cardinals) contain $\kappa$-standard morasses.
\smallskip\\
This is a part of my dissertation [Irr1]. I thank Dieter Donder for being my adviser, Hugh Woodin for an invitation to Berkeley, where part of the work was done, and the DFG-Graduiertenkolleg ``Sprache, Information, Logik'' in Munich for their support. 

\section{$(\omega _1 , \beta )$-Morasses}
Let me briefly recall how an object of size $\omega_2$ is constructed from countable objects in G\"odel's constructible universe $L$. That is, let me briefly describe how a construction along an $(\omega_1,1)$-morass works. Let an ordinal $\nu$ be called $\omega_2$-like, if $L_\nu \models$ there exists exactly one uncountable cardinal. Let $S^0=\{ \alpha \in Lim \mid L_\nu \models ( \alpha = \omega_1)$ for some $\omega_2$-like ordinal $\nu\}$. Then there are different kinds of $\omega_2$-like ordinals, namely for every $\alpha \in S^0$ there is the set $S_\alpha = \{ \nu\mid \nu$ is $\omega$-like and $L_\nu \models \alpha=\omega_1\}$ of  those which believe that $\alpha = \omega_1$. Now, a morass construction proceeds as follows: On the one hand, one constructs for every $\alpha \in S^0 \cap \omega_1$ by induction over $\nu \in S_\alpha$ a  countable chain $\langle \frak{A}_\nu \mid \nu \in S_\alpha\rangle$ of countable structures $\frak{A}_\nu$. On the other hand, one constructs by induction over $\alpha \in S^0$ a system of embeddings between these chains. As direct limit of this system of embeddings, one obtains a chain  $\langle \frak{A}_\nu \mid \nu \in S_{\omega_1}\rangle$ of length $\omega_2$ of structures $\frak{A}_\nu$ of size $\leq \omega_1$. Finally, the structure $\frak{A}$ of size $\omega_2$ that one wants to construct is obtained as the direct limit of this chain. 
\smallskip\\
The approach is generalized by Jensen to all $\beta < \omega_1$.  Let an ordinal $\nu$ be $\omega_{1+\beta}$-like, if the set $\{ \alpha \mid L_\nu \models \alpha$ is an uncountable cardinal$\} \cup \{ \nu\}$ has order-type $\beta+1$. The basic construction is first carried out for countable structures $\frak{A}_\nu$ and all $\omega_{1+\beta}$-like ordinals $\nu$ with $\nu < \omega_1$, and then directed systems of embeddings are used to blow it up to $\omega_{1+\beta}$. This motivates his definition of $(\omega_\alpha,\beta)$-morasses. They describe axiomatically the properties of the $\omega_{1+\beta}$-like ordinals which enable such constructions. This short introduction to morasses explains already why Jensen never introduced $(\omega_\alpha,\beta)$-morasses for the case $\omega_1 \leq \beta$, namely because then there are no $\omega_{1+\beta}$-like ordinals below $\omega_1$.  
 \smallskip\\
 To explain how I circumvent this problem, let me first introduce the notation $f:\bar\nu \Rightarrow \nu$ from Jensen's approach. As  explained, he considers the sets $S_\alpha =\{ \nu \mid L_\nu \models \alpha$ is the largest cardinal$\}$. Let $\alpha_\nu $ be the largest cardinal of $L_\nu$.  Then he constructs, on the one hand, by induction over the ordinals in the sets $S_\alpha$ chains  $\langle \frak{A}_\nu \mid \nu \in S_\alpha\rangle$ of structures $\frak{A}_\nu$. On the other hand, he considers maps $f$ which map under certain conditions $S_{\alpha_{\bar\nu}} \cap \bar \nu$ into $S_{\alpha_\nu}\cap \nu$ in such a way, that $f$ can be extended to an embedding from $\langle \frak{A}_\tau \mid \tau \in S_{\alpha_{\bar\nu}}\cap \bar \nu\rangle$ into $\langle \frak{A}_\tau \mid \tau \in S_{\alpha_\nu}\cap\nu\rangle$. For such a map he uses the notation $f: \bar \nu \Rightarrow \nu$. The possibility to extend the maps to uniform constructions is guaranteed by the so-called logical preservation axiom (see axiom LP1 below). If $f:\bar\nu\Rightarrow \nu$, then in Jensen's case $\bar\nu$ and $\nu$ are of the same type, that is, if $\bar\nu$ is $\omega_\eta$-like, then $\nu$ is also $\omega_\eta$-like. In my case, they can have different types, i.e. if $\bar\nu$ is $\omega_\eta$-like, then $\nu$ can be $\omega_\gamma$-like for some $\gamma \geq \eta$. This is done in such a way that in the limit a construction along the $\omega_{1+\beta}$-like cardinals takes place.
 \smallskip\\
 As consequence, also the chains $\langle \frak{A}_\tau \mid \tau \in S_{\alpha_{\bar\nu}}\cap \bar \nu\rangle$ and $\langle \frak{A}_\tau \mid \tau \in S_{\alpha_\nu}\cap\nu\rangle$ of ordinals $\tau$ of different types will have to fit together. The idea is to take care of this in the recursive definition of the chains. However, this makes it necessary to introduce a second logical preservation axiom which guarantees  that if $\bar \nu$ and $\nu$ are of different types and $f:\bar\nu \Rightarrow \nu$ is cofinal, then $f$ can be extended to an embedding from  $\langle \frak{A}_\tau \mid \tau \in S_{\alpha_{\bar\nu}}\cap \bar \nu\rangle$ into $\langle \frak{A}_\tau \mid \tau \in S_{\alpha_\nu}\cap\nu\rangle$. The second logical preservation axiom (see LP2 below) is inspired by and closely related to the construction of $\Box$-sequences. Unfortunately, I do not have an example of a typical recursive morass construction which can be carried out with my morasses but not with Jensen's morasses or a $\Box$-principle.
 \smallskip\\
 Morasses are also closely related to Jensen's fine structure theory in the following way. If for a map $f$ the relation $f:\bar\nu\Rightarrow\nu$ holds, this does not only mean that $f:\bar\nu\rightarrow\nu$ but also that it can be extended to a map $f:\mu_{\bar\nu}\rightarrow\mu_\nu$ where $\mu_{\bar\mu} \geq \bar \nu$ and $\mu_\nu \geq \nu$ depend on $\bar\nu$ and $\nu$. In this sense, it can be interpreted as saying not only that $f$ is a map from $\bar\nu$ to $\nu$, but that it is a $\Sigma_1$-elementary map from $\langle L_{\bar\nu},\bar A\rangle$ into $\langle L_\nu,A\rangle$ where $\bar A$ is a predicate coding $L_{\mu_{\bar\nu}}$ and $A$ is a predicate coding $L_{\mu_\nu}$. To show the above mentioned property that from my morasses an inner model with fine structure can be constructed, I include into my definition of morasses explicitly such a coding property. Moreover, the fine structural coding property for $\Sigma_n$-elementary maps is represented by relations $\Rightarrow_n$.
 \bigskip\\
 Let $\omega _1 \leq \beta$, $S=Lim \cap \omega _{1+\beta}$ and $\kappa :=\omega _{1+\beta}$. 
\smallskip\\
 We write $Card$ for the class of cardinals and $RCard$ for the class of regular cardinals.
\smallskip\\
Let $\vartriangleleft$ be a binary relation on $S$ such that:
\smallskip \\
(a) If $\nu \vartriangleleft \tau$, then $\nu < \tau$.

For all $\nu \in S - RCard$, $\{ \tau \mid \nu \vartriangleleft \tau \}$ is closed. 

For $\nu \in S - RCard$, there is a largest $\mu$ such that $\nu \trianglelefteq \mu$.
\smallskip \\
Let $\mu _\nu$ be this largest $\mu$ with $\nu \trianglelefteq \mu$.
\smallskip\\
Let 
$$\nu \sqsubseteq \tau : \Leftrightarrow \nu \in Lim(\{ \delta \mid \delta \vartriangleleft \tau \} ) \cup \{ \delta \mid \delta \trianglelefteq \tau \} .$$
(b) $\sqsubseteq$ is a (many-rooted) tree. 
\smallskip\\
Hence, if $\nu \notin RCard$ is a successor in $\sqsubset$, then $\mu _\nu$ is the largest $\mu$ such that $\nu \sqsubseteq \mu$. To see this, let $\mu ^\ast _\nu$ be the largest $\mu$ such that $\nu \sqsubseteq \mu$. It is clear that $\mu _\nu \leq \mu _\nu ^\ast$, since $\nu \trianglelefteq \mu$ implies $\nu \sqsubseteq \mu$. So assume that $\mu _\nu < \mu _\nu ^\ast$. Then $\nu \not\vartriangleleft \mu_\nu^\ast$ by the definition of $\mu _\nu$. Hence $\nu \in Lim(\{ \delta \mid \delta \vartriangleleft \mu _\nu ^\ast\} )$ and $\nu \in Lim(\{ \delta \mid \delta \sqsubseteq \mu _\nu^\ast\} )$. Therefore, $\nu \in Lim(\sqsubseteq)$ since $\sqsubseteq$ is a tree. That contradicts our assumption that $\nu$ is a successor in $\sqsubset$.
\smallskip\\
The properties of $\omega\nu \vartriangleleft \omega\tau$ are an axiomatic description of the relation "$\omega\nu$ is regular in $L_\tau$". If $\omega\nu \vartriangleleft \omega\tau$ really is this relation, then $\omega\nu \sqsubset \omega\tau$ implies that $\omega\nu$ is a cardinal in $L_\tau$, while the converse implication is not true in general. This is a crucial difference to Jensen's morasses, where $\omega\nu \sqsubset \omega\tau$ is an axiomatic description of "$\omega\nu$ is a cardinal in $L_\tau$", and it is the reason why $\vartriangleleft$ is introduced. However, if there exists a maximal cardinal in $L_\nu$ and $\nu < \tau$, then the two interpretations of $\omega \nu \sqsubset \omega\tau$ coincide. 
\smallskip \\
For $\alpha \in S$, let $|\alpha|$ be the rank of $\alpha$ in this tree. Let 
\smallskip 

$S^+ := \{ \nu \in S \mid \nu$ is a successor in $\sqsubset \}$
\smallskip

$S^0 := \{ \alpha \in S \mid |\alpha| =0\}$
\smallskip

$\widehat{S^+}:=\{ \mu _\tau \mid \tau \in S^+ - RCard \}$
\smallskip

$\widehat{S}:=\{ \mu _\tau \mid \tau \in S -RCard \}$.
\smallskip \\
Let $S_\alpha := \{ \nu \in S \mid \nu$ is a direct successor of $\alpha$ in $\sqsubset \}$. For $\nu \in S^+$, let $\alpha _\nu$ be the direct predecessor of $\nu$ in $\sqsubset$. For $\nu \in S^0$, let $\alpha _\nu :=0$. For $\nu \not\in S^+ \cup S^0$, let $\alpha _\nu :=\nu$.
\smallskip \\
(c) For $\nu ,\tau \in (S^+ \cup S^0)-RCard$ such that $\alpha _\nu =\alpha _\tau$, suppose:
$$\nu < \tau \quad \Rightarrow \quad \mu _\nu < \tau .$$
For all $\alpha \in S$, suppose:
\begin{tabbing}
(d) \= $S_\alpha$ is closed \\[0.25ex]
(e) \> $card(S_\alpha ) \leq \alpha ^+$ \\[0.25ex]
\> $card(S _\alpha ) \leq card(\alpha )$ if $card(\alpha ) < \alpha$ \\[0.25ex]
(f) \> $\omega _1 = max(S^0 ) = sup(S^0 \cap \omega _1 )$ \\[0.25ex]
\> $\omega _{1+i+1}=max(S_{\omega _{1+i}})=sup(S_{\omega _{1+i}} \cap \omega _{1+i+1})$ for all $i<\beta$. \\
\end{tabbing}
Let $D=\langle D_\nu \mid \nu \in \widehat{S} \rangle$ be a sequence such that $D_\nu \subseteq J^D_\nu$. To simplify matters, my definition of $J^D_\nu$ is such that $J^D_\nu \cap On =\nu$ (see section 3 or [SchZe]).
\medskip \\
Let an $\langle S,\vartriangleleft ,D\rangle$-maplet $f$ be a triple $\langle \bar \nu ,|f| , \nu \rangle$ such that $\bar \nu ,\nu \in S - RCard$ and $|f| :J^D_{\mu _{\bar \nu}} \rightarrow J^D_{\mu _\nu}$.
\smallskip \\
Let $f=\langle \bar \nu ,|f| , \nu \rangle$ be an  $\langle S,\vartriangleleft ,D\rangle$-maplet. Then we define $d(f)$ and $r(f)$ by $d(f)=\bar \nu$ and $r(f)=\nu$. Set $f(x):=|f| (x)$ for $x \in J^D_{\mu _{\bar \nu}}$ and $f(\mu_{\bar\nu} ):=\mu _\nu$. But $dom(f)$, $rng(f)$, $f \upharpoonright X$, etc. keep their usual set-theoretical meaning, i.e. $dom(f)=dom(|f| )$, $rng(f)=rng(|f| )$, $f \upharpoonright X = |f| \upharpoonright X$, etc. 
\smallskip \\
For $\bar \tau \leq \mu_{ \bar \nu}$, let $f^{( \bar \tau )}=\langle \bar \tau , |f| \upharpoonright J^D_{\mu _{\bar \tau}}, \tau \rangle$ where $\tau =f( \bar \tau )$. Of course, $f^{( \bar \tau )}$ needs not to be a maplet. The same is true for the following definitions.  Let $f^{-1}=\langle \nu , |f| ^{-1}, \bar \nu \rangle$. For $g= \langle \nu , |g| , \nu ^\prime \rangle$ and $f=\langle \bar \nu ,|f| , \nu \rangle$, let $g \circ f= \langle \bar \nu , | g| \circ |f| , \nu ^\prime \rangle$. If $g= \langle \nu ^\prime , | g | , \nu \rangle$ and $f=\langle \bar \nu ,| f| , \nu \rangle$ such that $rng(f) \subseteq rng(g)$, then set $g^{-1}f=\langle \bar \nu , | g| ^{-1} \mid f \mid , \nu ^\prime \rangle$. Finally set $id _\nu = \langle \nu , id \upharpoonright J^D_{\mu _\nu} , \nu \rangle$.
\smallskip \\
Let $\frak{F}$ be a set of  $\langle S,\vartriangleleft ,D\rangle$-maplets $f=\langle \bar \nu ,|f| , \nu \rangle$ such that the following \\
holds:
\smallskip \\
(0) $f(\bar \nu )=\nu$, $f(\alpha _{\bar \nu})=\alpha _\nu$ and $|f|$ is order-preserving. 
\smallskip \\
(1) For $f \neq id _{\bar \nu}$, there is some $\beta \sqsubseteq \alpha _{\bar \nu}$ such that $f \upharpoonright \beta = id \upharpoonright \beta$ and $f(\beta )>\beta$.
\smallskip \\
(2) If $\bar \tau \in S^+$ and $\bar \nu \sqsubset \bar \tau \sqsubseteq \mu_{\bar \nu}$, then $f^{(\bar \tau )}\in \frak{F}$.
\smallskip \\
(3) If $f,g \in \frak{F}$ and $d(g)=r(f)$, then $g \circ f \in \frak{F}$.
\smallskip\\
(4) If $f,g \in \frak{F}$, $r(g)=r(f)$ and $rng(f)\subseteq rng(g)$, then $g^{-1} \circ f \in \frak{F}$.
\smallskip\\
We write $f: \bar \nu \Rightarrow \nu$ if $f=\langle \bar \nu , |f| , \nu \rangle \in \frak{F}$. If $f \in \frak{F}$ and $r(f)=\nu$, then we write $f \Rightarrow \nu$. The uniquely determined $\beta$ in (1) shall be denoted by $\beta (f)$.
\smallskip \\
Say $f \in \frak{F}$ is minimal for a property $P(f)$ if $P(f)$ holds and $P(g)$ implies $g^{-1}f \in \frak{F}$. 
\smallskip \\
Let 
\smallskip

$f_{(u,x ,\nu )}=$ the unique minimal $f \in \frak{F}$ for $f \Rightarrow \nu$ and $u \cup \{ x \} \subseteq rng(f)$,
\smallskip \\
if such an $f$ exists. The axioms of the morass will guarantee that $f_{(u,x ,\nu )}$ always exists if $\nu \in S-RCard^{L_\kappa[D]}$. Therefore, we will always assume and explicitly mention that $\nu \in S-RCard^{L_\kappa[D]}$ when $f_{(u,x ,\nu )}$ is mentioned.
\medskip \\
Say $\nu \in S-RCard^{L_\kappa[D]}$ is independent if $d(f_{(\beta ,0,\nu )}) < \alpha _\nu$ holds for all $\beta < \alpha _\nu$.
\smallskip \\
For $\tau \sqsubseteq \nu \in S-RCard^{L_\kappa[D]}$, say $\nu$ is $\xi$-dependent on $\tau$ if $f_{(\alpha _\tau , \xi , \nu )} = id _\nu$.
\medskip \\
For $f \in \frak{F}$, let $\lambda (f):= sup(f[d(f)])$.
\medskip\\
For $\nu \in S-RCard^{L_\kappa[D]}$, let
$$C_\nu =\{ \lambda (f) <\nu \mid f \Rightarrow \nu \} $$
$$\Lambda (x ,\nu )=\{ \lambda (f_{(\beta , x , \nu )}) < \nu \mid \beta < \nu \} .$$
It will be shown that $C_\nu$ and $\Lambda (x ,\nu )$ are closed in $\nu$.
\smallskip\\
Recursively define a function $q_\nu:k_\nu +1 \rightarrow On$, where $k_\nu \in \omega$:
\smallskip 

$q_\nu (0)=0$
\smallskip

$q_\nu (k+1) = max(\Lambda (q_\nu \upharpoonright (k+1) ,\nu ))$
\smallskip \\
if $max(\Lambda (q_\nu \upharpoonright (k+1) ,\nu ))$ exists. The axioms will guarantee that this recursion breaks off (see lemma 4 below), i.e. there is some $k_\nu$ such that either
\smallskip

$\Lambda (q_\nu \upharpoonright ({k_\nu}+1),\nu ) =\emptyset$
\smallskip \\
or
\smallskip

$\Lambda (q_\nu \upharpoonright ({k_\nu}+1),\nu )$ is unbounded in $\nu$.
\medskip \\
Define by recursion on $1 \leq n \in \omega$, simultaneouly for all $\nu \in S-RCard^{L_\kappa[D]}$,  $\beta \in \nu$ and $x \in J^D_{\mu_\nu}$ the following notions. Here definitions are to be understood in Kleene's sense, i.e., that the left side is defined iff the right side is, and in that case, both are equal.
\smallskip

$f^1_{(\beta ,x,\nu )}=f_{(\beta ,x ,\nu )}$
\smallskip

$\tau (n,\nu )$ $=$ the least $\tau \in S^0 \cup S^+ \cup \widehat{S}$ such that for some $x \in J^D_{\mu _\nu}$
$$f^n_{(\alpha _\tau ,x ,\nu )}=id_\nu$$

$x (n,\nu )$ $=$ the least $x \in J^D_{\mu _\nu}$ such that $f^n_{(\alpha _{\tau (n,\nu )} ,x ,\nu )}=id_\nu$
\smallskip 

$K^n_\nu =\{ d(f^n_{(\beta ,x (n,\nu ),\nu )})< \alpha _{\tau (n,\nu )} \mid \beta < \nu \}$
\smallskip

$f \Rightarrow _n \nu$ iff $f \Rightarrow \nu$ and for all $1 \leq m < n$ $$rng(f) \cap J^D_{\alpha _{\tau (m,\nu )}} \prec _1 \langle J^D_{\alpha _{\tau (m,\nu )}},D \upharpoonright {\alpha _{\tau (m,\nu )}},K^m_\nu \rangle $$
$$x (m,\nu ) \in rng(f)$$

$f^n_{(u ,\nu )}$ $=$ the minimal $f \Rightarrow _n \nu$ such that $u \subseteq rng(f)$
\smallskip 

$f^n_{(\beta , x ,\nu )}$ $=$ $f^n_{(\beta \cup \{ x \} ,\nu )}$
\smallskip 

$f:\bar \nu \Rightarrow _n \nu$ $:\Leftrightarrow$ $f \Rightarrow _n \nu$ $and$ $f:\bar \nu \Rightarrow \nu$.
\medskip\\
Let
\smallskip

$n_\nu$ $=$ the least $n$ such that $f^n_{(\gamma ,x ,\mu _\nu )}$ is confinal in $\nu$ for some $x \in J^D_{\mu _\nu}$, $\gamma \sqsubset \nu$
\smallskip

$x _\nu$ $=$ the least $x$ such that $f^{n_\nu}_{(\alpha _\nu ,x ,\mu _\nu )}=id_{\mu _\nu}$.
\smallskip \\
Let
\smallskip 

$\alpha ^\ast _\nu =\alpha _\nu$ if $\nu \in S^+$
\smallskip 

$\alpha _\nu ^\ast = sup \{ \alpha < \nu \mid \beta (f^{n_\nu}_{(\alpha ,x _\nu ,\mu _\nu )})=\alpha \}$ if $\nu \notin S^+$.
\smallskip\\
Let $P_\nu := \{ x _\tau \mid \nu \sqsubset \tau \sqsubseteq \mu _\nu ,\tau \in S^+ \} \cup \{ x _\nu \}$.
\bigskip \\
We say that $\frak{M} = \langle S,\vartriangleleft ,\frak{F},D \rangle$ is an $(\omega _1 , \beta )$-morass if the following axioms hold:
\medskip \\
{\bf (MP -- minimum principle)}
\smallskip\\
If $\nu \in S-RCard^{L_\kappa [D]}$ and $x \in J^D_{\mu _\nu}$, then $f_{(0,x ,\nu)}$ exists.
\medskip \\
{\bf (LP1 -- first logical preservation axiom)}
\smallskip\\
If $f:\bar \nu \Rightarrow \nu$, then $|f| : \langle J^D_{\mu _{\bar \nu}} , D \upharpoonright {\mu _{\bar \nu}} \rangle \rightarrow \langle J^D_{\mu _\nu} , D \upharpoonright {\mu _\nu} \rangle$ is $\Sigma _1$-elementary.
\medskip\\
{\bf (LP2 -- second logical preservation axiom)}
\smallskip\\
Let $f:\bar \nu \Rightarrow \nu$ and $f(\bar x )=x$. Then $$(f \upharpoonright J^D_{\bar \nu }):\langle J^D_{\bar \nu} ,D \upharpoonright {\bar \nu},\Lambda (\bar x , \bar \nu ) \rangle \rightarrow \langle J^D_\nu ,D \upharpoonright \nu ,\Lambda (x ,\nu ) \rangle$$ is $\Sigma _0$-elementary.
\medskip \\
{\bf (CP1 -- first continuity principle)}
\smallskip\\
For $i \leq j < \lambda$, let $f_i:\nu _i \Rightarrow \nu$ and $g_{ij}:\nu _i \Rightarrow \nu _j$ such that $g_{ij}=f^{-1}_jf_i$. Let $\langle g_i \mid i<\lambda \rangle$ be the transitive, direct limit of the directed system $\langle g_{ij} \mid i \leq j<\lambda \rangle$ and $hg_i= f_i$ for all $i<\lambda$. Then $g_i,h \in \frak{F}$.
\medskip \\
{\bf (CP2 -- second continuity principle)}
\smallskip\\
Let $f:\bar \nu \Rightarrow \nu$ and $\lambda =sup(f[\bar \nu ])$. If, for some $\bar \lambda$, $h:\langle J^{\bar D}_{\bar \lambda} , \bar D \rangle \rightarrow \langle J^D_\lambda, D \upharpoonright \lambda \rangle$ is $\Sigma _1$-elementary and $rng(f \upharpoonright J^D_{\bar \nu} ) \subseteq rng(h)$, then there is some $g:\bar \lambda \Rightarrow \lambda$ such that $g \upharpoonright J^{\bar D}_{\bar \lambda} = h$. 
\medskip \\
{\bf (CP3 -- third continuity principle)}
\smallskip \\
If $C_\nu =\{ \lambda (f) < \nu \mid f \Rightarrow \nu \}$ is unbounded in $\nu \in S-RCard^{L_\kappa[D]}$,  then the following holds for all $x \in J^D_{\mu _\nu}$:
$$rng(f_{(0,x ,\nu )})=\bigcup \{ rng(f_{(0,x ,\lambda )}) \mid \lambda \in C_\nu \} .$$
{\bf (DP1 -- first dependency axiom)}
\smallskip\\
If $\mu _\nu < \mu _{\alpha _\nu}$, then $\nu \in S-RCard^{L_\kappa[D]}$ is independent.
\medskip \\
{\bf (DP2 -- second dependency axiom)}
\smallskip \\
If $\nu\in S-RCard^{L_\kappa[D]}$ is $\eta$-dependent on $\tau \sqsubseteq \nu$, $\tau \in S^+$,  $f:\bar \nu \Rightarrow \nu$, $f(\bar \tau )=\tau$ and $\eta \in rng(f)$, then $f^{(\bar \tau )}: \bar \tau \Rightarrow \tau$.
\medskip \\
{\bf (DP3 -- third dependency axiom)}
\smallskip\\
If $\nu \in \widehat{S}-RCard^{L_\kappa[D]}$ and $1 \leq n \in \omega$, then the following holds:
\smallskip \\
(a) If $f^n_{(\alpha _\tau , x ,\nu )}=id_\nu$, $\tau \in S^+ \cup S^0$ and $\tau \sqsubseteq \nu$, then $\mu _\nu =\mu _\tau$.
\smallskip \\
(b) If $\beta < \alpha _{\tau (n,\nu )}$, then also $d(f^n_{(\beta ,x (n,\nu ),\nu )})< \alpha _{\tau (n,\nu )}$.
\medskip \\
{\bf (DF -- definability axiom)}
\smallskip \\
(a) If $f_{(0,z_0 ,\nu )}=id_\nu$ for some $\nu \in \widehat{S}-RCard^{L_\kappa[D]}$ and $z_0 \in J^D_{\mu _\nu}$, then
$$\{ \langle z ,x,f_{(0,z ,\nu )}(x) \rangle \mid z \in J^D_{\mu _\nu}, x \in dom(f_{(0,z ,\nu )}) \}$$
is uniformly  definable over $\langle J^D_{\mu _\nu},D\upharpoonright {\mu _\nu},D_{\mu _\nu} \rangle$.
\smallskip \\
(b) For all $\nu \in S-RCard^{L_\kappa [D]}$, $x \in J^D_{\mu _\nu}$, the following holds:
$$f_{(0,x ,\nu )}=f^{n_\nu}_{(0,\langle x ,\nu , \alpha ^\ast _\nu ,P_\nu \rangle ,\mu _\nu )}.$$
\smallskip \\
This finishes the definition of an $(\omega_1,\beta)$-morass.
\bigskip\\
A consequence of the axioms is ($\times$):
\pagebreak\\
{\bf Theorem}
\smallskip\\
$$\{ \langle z ,\tau ,x,f_{(0,z ,\tau )}(x) \rangle \mid  \tau < \nu ,\mu _\tau =\nu , z \in J^D_{\mu _\tau}, x \in dom(f_{(0,z ,\tau )}) \}$$
$$\cup \{ \langle z ,x,f_{(0,z ,\nu )}(x) \rangle \mid \mu _\nu =\nu ,z \in J^D_{\mu _\nu}, x \in dom(f_{(0,z ,\nu )}) \}$$
$$\cup (\sqsubset \cap \nu ^2)$$
is for all $\nu \in S$ uniformly definable over $\langle J^D_\nu ,D\upharpoonright \nu ,D_\nu \rangle$.
\bigskip\\
The proof of the property ($\times$) streches over the next twelve lemmas unto the end of the section. It is proved by induction over $\mu \in \widehat{S}$, i.e. we prove it for all $\nu$ with $\mu_\nu=\mu$ assuming that it holds for all $\tau$ with $\mu_\tau<\mu$. More precisely, assume it holds for all $\tau$ such that $\mu _\tau < \mu$. Then we show that the various minimal maps $f^n_{(u,\nu )}$ exist for all $\nu$ such that $\mu _\nu =\mu$ and all $u \subseteq \mu$ (lemmas 1 and 12). And we show that $q_\nu$ exists for all $\nu$ such that $\mu _\nu =\mu$ (lemma 4). Finally we prove that ($\times$) holds for all $\nu$ with $\mu _\nu =\mu$.
\medskip\\
So assume ($\times$) for all $\tau$ such that $\mu _\tau < \mu _\nu := \mu$. If $\mu=0$, this holds trivially, because then there are no such $\tau$. For the proof we need the following lemmas which are very important in themselves but proved as part of our big induction on $\mu _\nu$.
\medskip\\ 
{\bf Lemma 1}
\smallskip \\
Let $\nu \in S-RCard^{L_\kappa [D]}$ and $u \subseteq J^D_{\mu _\nu}$. Then there is a minimal $f \in \frak{F}$ for $f \Rightarrow \nu$ and $u \subseteq rng(f)$.
\smallskip \\
We write $f_{(u,\nu)}$ for this $f$.
\smallskip\\
{\bf Proof:}
\smallskip \\
(1) For finite $u=\{ \xi _1 ,\dots ,\xi _n \}$, we have $f_{(u,\nu )}=f_{(0,\langle \xi _1, \dots ,\xi _n \rangle ,\nu )}$.
\smallskip \\
For, by (LP1), $f_{(u,\nu )}:\langle J^D_{\mu _{\bar \nu _1}},D \upharpoonright \mu _{\bar \nu _1} \rangle \rightarrow \langle J^D_{\mu_\nu} ,D \upharpoonright \nu \rangle$ is $\Sigma _1$-elementary. Since $J^D_{\mu_\nu}$ is closed under pairs, $u \subseteq rng(f_{(u,\nu )})$ implies $\langle \xi _1, \dots ,\xi _n \rangle \in rng(f_{(u,\nu )})$. For the converse, we note that $f_{(0,\langle \xi _1, \dots ,\xi _n \rangle ,\nu )}:\langle J^D_{\mu_{\bar \nu _2}},D\upharpoonright \mu_{\bar \nu _2} \rangle \rightarrow \langle J^D_{\mu_\nu} ,D\upharpoonright \mu_\nu \rangle$ is  $\Sigma _1$-elementary by (LP1). Hence $\langle \xi _1, \dots ,\xi _n \rangle \in rng(f_{(0,\langle \xi _1, \dots ,\xi _n \rangle ,\nu )})$ implies $u \subseteq rng(f_{(0,\langle \xi _1, \dots ,\xi _n \rangle ,\nu )})$. By (MP), $f_{(0,\{ \xi _1 ,\dots ,\xi _n \},\nu)}$ exists, and by its minimality, it is as wished.
\smallskip \\
(2) Now, let $u$ be infinite. Then $I=\{ v \subseteq u \mid v$ finite $\}$ is directed with regard to $\subseteq$. Let $g_{vw}=f^{-1}_{(w,\nu )} f_{(v,\nu )}$ for $v \subseteq w \in I$. Then $g_{vw} \in \frak{F}$ by (4) and the definition of minimality. Let $\langle g_v \mid v \in I \rangle$ be the transitive, direct limit of $\langle g_{vw} \mid v \subseteq w \rangle$ and $hg_v= f_{(v,\nu)}$ for all $v \in I$. Then $g_v,h \in \frak{F}$ by (CP1). But obviously $h=f_{(u,\nu )}$. $\Box$
\medskip \\
{\bf Lemma 2}
\smallskip \\
Let $\nu \in S-RCard^{L_\kappa [D]}$. Then:
\smallskip \\
(a) Let $g:\bar \nu \Rightarrow \nu$, $\bar u \subseteq J^D_{\mu _{\bar\nu}}$ and $u=g[\bar u]$. Then $gf_{(\bar u ,\bar \nu )}=f_{(u,\nu )}$.
\smallskip \\
(b) $id_\nu \in \frak{F}$.
\smallskip\\
(c) If $f \Rightarrow \nu$ and $f \upharpoonright \alpha _\nu =id \upharpoonright \alpha _\nu$, then $f=id_\nu$.
\smallskip\\
(d) $J^D_{\mu _\nu} = \bigcup \{ rng(f_{(\beta , \xi , \nu )}) \mid \beta < \alpha _\nu \}$ for all $\xi \in  J^D_{\mu _\nu}$.
\smallskip \\
{\bf Proof:}
\smallskip \\
(a) On the one hand, we have
\smallskip 

$\bar u = g^{-1}[u] \subseteq rng(g^{-1}f_{(u,\nu )})$
\smallskip

$\Rightarrow$ $rng(f_{(\bar u,\bar \nu )}) \subseteq rng(g^{-1}f_{(u,\nu )})$
\smallskip

$\Rightarrow$ $rng(gf_{(\bar u,\bar \nu )}) \subseteq rng(f_{(u,\nu )})$.
\smallskip \\
On the other hand, we have
\smallskip 

$u \subseteq rng(gf_{(\bar u,\bar \nu )})$
\smallskip 

$\Rightarrow$ $rng(f_{(u,\nu )}) \subseteq rng(gf_{(\bar u,\bar \nu )})$.
\smallskip \\
(b) $id_\nu =f_{(u ,\nu )}$ where $u=J^D_{\mu _\nu}$.
\smallskip \\
(c) Assume $f \neq id_\nu$. Then $\beta (f)\leq \alpha _\nu$ by axiom (1). But $f \upharpoonright \alpha _\nu =id \upharpoonright \alpha _\nu$ by the hypothesis and $f(\alpha _\nu )=\alpha _\nu$ by axiom (0). Contradiction!
\smallskip\\
(d) If we let $h:\bar \nu \Rightarrow \nu$ be the uncollapse of $\bigcup \{ rng(f_{(\beta , \xi , \nu )}) \mid \beta < \alpha _\nu \}$, then $h \in \frak{F}$ and $h \upharpoonright \alpha _\nu = id \upharpoonright \alpha _\nu$. So $h=id_\nu$ by (c). $\Box$
\medskip \\
{\bf Lemma 3}
\smallskip \\
Let $\bar \nu ,\nu \in S$ and let $h:\langle J^{\bar D}_{\bar \nu} , \bar D \rangle \rightarrow \langle J^D_\nu ,D\upharpoonright\nu \rangle$ be $\Sigma _1$-elementary such that there is some $\beta \sqsubseteq \bar \nu$ with $h \upharpoonright \beta = id \upharpoonright \beta$. Let $h(\mu _{\bar \tau})=\mu _\tau < \nu$ and $\tau =h(\bar \tau ) \in S-RCard^{L_\kappa [D]}$. Then $h^{(\bar \tau )}: \bar \tau \Rightarrow \tau$.
\smallskip \\
{\bf Proof:} Let $\delta _{\bar\tau} \sqsubseteq \bar\tau$ and $\delta _{\bar\nu} \sqsubseteq \bar\nu$ be minimal. If $\delta _{\bar\tau} \not\sqsubseteq \bar\nu$, then $\mu _{\bar\tau} < \delta _{\bar\nu}$. To see this, we consider the three cases $\delta _{\bar\tau} = \delta _{\bar\nu}$, $\delta _{\bar\tau} > \delta _{\bar\nu}$ and $\delta _{\bar\tau} < \delta _{\bar\nu}$. The first case is impossible because if $\delta _{\bar\nu} =\delta _{\bar\tau}$, then $\delta _{\bar\tau} \sqsubseteq \bar\nu$ by definition of $\delta _{\bar\nu}$. The second case is impossible because then by axiom (c) $\mu _{\delta _{\bar\nu}} < \delta _{\bar\tau}$. But $\delta _{\bar\tau} \leq \bar\tau$ by definition of $\delta _{\bar\tau}$ and $\mu _{\bar\nu} \leq \mu _{\delta _{\bar\nu}}$ by definition of $\delta _{\bar\nu}$ and $\mu _\cdot$. Hence $\mu _{\bar\nu} \leq \mu _{\delta _{\bar\nu}} < \delta _{\bar\tau} \leq \bar\tau$ which contradicts the assumption $\mu _{\bar\tau} < \bar\nu$. Hence $\delta _{\bar\tau} < \delta _{\bar\nu}$ must hold. But then $\mu _{\delta _{\bar\tau}} < \delta _{\bar\nu}$ by axiom (c) and therefore $\mu _{\bar\tau} \leq \mu _{\delta _{\bar\tau}} < \delta _{\bar\nu}$ as claimed, by definition of $\delta _{\bar\tau}$ and $\mu _\cdot$. So by assumption $h^{(\bar\tau)}=id_{\bar\tau}$ and $id_{\bar\tau} \in \frak{F}$ by lemma 2 (b).
\smallskip\\
Now, let $\bar\delta := \delta _{\bar\tau} \sqsubseteq \bar\nu$ and $f_{(\bar\delta ,x ,\bar\tau )}:\bar\tau (x ) \Rightarrow \bar\tau$. Let $\bar\delta \sqsubset \bar\gamma (x ) \sqsubseteq \bar\tau (x )$ where $\alpha _{\bar\gamma (x )} =\bar\delta$. Then, by (DP2), $f_{(0,x ,\bar\tau (x ) )}=f_{(0,x ,\bar\gamma (x ))}$ for all $x \in J^D_{\mu _{\bar\tau}}$. And we get $\mu _{\bar\gamma (x )} \leq \mu  _{\bar\tau} < \bar\nu \leq \mu _{\bar\delta}$. So, by (DP1), $\bar\gamma (x )$ is independent. That is, $d(f_{(\beta ,0,\bar\gamma (x ))}) < \alpha _{\bar\gamma (x )}$ for all $\beta < \alpha _{\bar\gamma (x )}$. Since $J^D_{\mu _{\bar\gamma (x )}}=\bigcup \{ rng(f_{(\beta ,0,\bar\gamma (x ))}) \mid \beta < \alpha _{\bar\gamma (x )}\}$, $x \in rng(f_{(\beta ,0,\bar\gamma (x ))})$ for some $\beta < \alpha _{\bar\gamma (x )}$. Hence $d(f_{(0,x ,\bar\gamma (x ))}) < \bar\delta$. Altogether, we get
\smallskip 

$d(f_{(0,x ,\bar\tau )})=d(f_{(\bar\delta ,x ,\bar\tau )}\circ f_{(0,x ,\bar\tau (x ))})=d(f_{(0,x ,\bar\tau (x ))})=d(f_{(0,x ,\bar\gamma (x ))})< \bar\delta$.
\smallskip \\
By our assuption $h \upharpoonright \bar \delta = id \upharpoonright \bar\delta$. And by our induction hypothesis, $(\times)$ holds for $\mu _\tau$. So by the $\Sigma _1$-elementarity of $h:\langle J^D_{\bar \nu} , \bar D \rangle \rightarrow \langle J^D_\nu ,D\upharpoonright \nu \rangle$, if $x \in rng(h)$, then even $rng(f_{(0,x ,\tau )}) \subseteq rng(h)$. Thus 
\smallskip 

$rng(h) \cap J^D_{\mu _\tau}= \bigcup \{ rng(f_{(0,x ,\tau )}) \mid x \in rng(h) \cap J^D_{\mu _\tau} \}$.
\smallskip \\
Therefore,
\smallskip 

$h^{(\bar \tau )}=f_{( u,\tau )} \in \frak{F}$ where $u=rng(h)\cap J^D_{\mu _\tau}$. $\Box$ 
\medskip\\
{\bf Lemma 4}
\smallskip\\
For all $\nu \in S-RCard^{L_\kappa [D]}$, $q_\nu$ exists.
\smallskip\\
{\bf Proof:} Suppose $q_\nu (k+1)=max (\Lambda (q_\nu \upharpoonright (k+1) ,\nu ))$ exists. Then $q_\nu (k+1) \in \Lambda (q_\nu \upharpoonright (k+1) ,\nu )$ and there is some $\beta$ such that $\lambda (f_{(\beta ,q_\nu \upharpoonright (k+1) ,\nu )})=q_\nu (k+1)$. The set of such $\beta$ is closed by (CP1). Thus there is a largest such $\beta$. Call it $\beta _k$. The recursion breaks off if the sequence $\langle \beta _k \mid k \rangle$ is strictly descending, since there is no descending sequence of length $\omega$. But $\beta _k \in rng(f_{(\beta _k,q_\nu \upharpoonright (k+2) ,\nu )})$ by ($\times$) and (LP2). Hence $\lambda (f_{(\beta _k,q_\nu \upharpoonright (k+2) ,\nu )})=\nu$ by the definition of $\beta _k$. Therefore, $\beta _{k+1}<\beta _k$. $\Box$
\medskip\\
{\bf Lemma 5}
\smallskip \\
Let $f:\bar \nu \Rightarrow \nu$, $x \in rng(f)$ and $\lambda =\lambda (f)$. Then $\Lambda (x ,\nu ) \cap \lambda =\Lambda (x ,\lambda )$.
\smallskip \\
{\bf Proof:} Let $f(\bar x )=x$. Then on the one hand, $(f \upharpoonright J^D_{\bar \nu}) : \langle J^D_{\bar \nu} ,D \upharpoonright {\bar \nu} ,\Lambda (\bar x ,\bar \nu ) \rangle \rightarrow \langle J^D_\nu ,D \upharpoonright \nu ,\Lambda (x ,\nu )\rangle$ is  $\Sigma _0$-elementary by (LP2). But then 
\smallskip 

$(\ast )$ $(f \upharpoonright J^D_{\bar \nu}) : \langle J^D_{\bar \nu} ,D \upharpoonright {\bar \nu} ,\Lambda (\bar x ,\bar \nu ) \rangle \rightarrow \langle J^D_\lambda ,D \upharpoonright \lambda ,\Lambda (x ,\nu ) \cap \lambda \rangle$ is also $\Sigma _0$-elementary.
\smallskip \\
On the other hand, by (CP2) and (LP2),
\smallskip 

$(\ast \ast )$ $(f \upharpoonright J^D_{\bar \nu}) : \langle J^D_{\bar \nu} ,D\upharpoonright {\bar \nu} ,\Lambda (\bar x ,\bar \nu ) \rangle \rightarrow \langle J^D_\lambda ,D\upharpoonright \lambda ,\Lambda (x ,\lambda ) \rangle$ is also $\Sigma _0$-elementary.
\smallskip \\
Consider the following three cases:
\smallskip \\
(1) $\Lambda (\bar x ,\bar \nu )=\emptyset$
\smallskip \\
Then, by $(\ast )$, $\Lambda (x ,\nu ) \cap \lambda =\emptyset$ and, by $(\ast \ast )$, $\Lambda (x ,\lambda )=\emptyset$.
\smallskip \\
(2) $\bar \eta := max(\Lambda (\bar x ,\bar \nu ))$ exists
\smallskip \\
Let $f(\bar \eta )=\eta$. Then, by $(\ast )$ and $(\ast \ast )$, 
\smallskip

$\eta =max(\Lambda (x ,\nu ) \cap \lambda )=max(\Lambda (x ,\lambda ))$.
\smallskip \\
And by (CP2), we have
\smallskip

$z \in \Lambda (\bar x ,\bar \nu )$ $\Leftrightarrow$ $z \in \Lambda (\bar x ,\bar \eta ) \cup \{ \bar \eta \}$.
\smallskip \\
But then, by $(\ast )$,
\smallskip

$z \in \Lambda (x ,\nu ) \cap \lambda$ $\Leftrightarrow$ $z \in \Lambda (x ,\eta ) \cup \{ \eta \}$.
\smallskip \\
and, because of $(\ast \ast )$,
\smallskip

$z \in \Lambda (x ,\lambda )$ $\Leftrightarrow$ $z \in \Lambda (x ,\eta ) \cup \{ \eta \}$. 
\smallskip \\
That's it!
\smallskip \\
(3) $\Lambda (\bar x ,\bar \nu )$ is unbounded in $\bar \nu$ 
\smallskip \\
Then, by $(\ast )$, $\Lambda (x ,\nu ) \cap \lambda$ is unbounded in $\lambda$. Hence $\lambda \in \Lambda (x ,\nu )$ because $\Lambda (x ,\nu )$ is closed. Therefore $\Lambda (x ,\lambda )=\Lambda (x ,\nu )\cap \lambda$ by (CP2). $\Box$
\bigskip \\
{\bf Lemma 6}
\smallskip \\
Let $f:\bar\nu \Rightarrow \nu$.
\smallskip \\
(a) If $q_\nu \upharpoonright k \in rng(f)$, then $f(q_{\bar\nu} \upharpoonright k)=q_\nu \upharpoonright k$.
\smallskip \\
(b) If $f$ is cofinal, then $f(q_{\bar \nu})=q_\nu$.
\smallskip \\
{\bf Proof:}
\smallskip \\
(a) That is proved by induction on $k$ using (LP2) to show $f(max(\Lambda (\bar x ,\bar \nu ))=max(\Lambda (x, \nu ))$ whenever $max(\Lambda (x ,\nu )) \in rng(f)$.
\smallskip \\
(b) Like (a). Since $f$ is cofinal, $q_\nu \upharpoonright (k+1)$ lies always in $rng(f)$. $\Box$ 
\medskip \\
{\bf Lemma 7}
\smallskip \\
$\lambda \in C_\nu$ implies $\lambda \in \Lambda (q_\lambda ,\nu )$.
\smallskip \\
{\bf Proof:} Since $\lambda \in C_\nu$, $q_\lambda \in rng(f)$ for some $f:\bar \nu \Rightarrow \nu$ by lemma 6 (b). So $\Lambda (q_\lambda ,\nu ) \cap \lambda =\Lambda (q_\lambda ,\lambda )$ by lemma 5. Therefore, by the definition of $q_\lambda$, $max(\Lambda (q_\lambda ,\nu ) \cap \lambda )$ does not exist. But if $\Lambda (q_\lambda ,\nu ) \cap \lambda$ is unbounded in $\lambda$, then $\lambda \in \Lambda (q_\lambda ,\nu )$ by the closedness of $\Lambda (q_\lambda ,\nu )$. So let $\Lambda (q_\lambda ,\nu ) \cap \lambda = \emptyset$. But then $\lambda =\lambda (f_{(0,q_\lambda ,\nu )})$. For $\lambda (f_{(0,q_\lambda ,\nu )}) \geq \lambda$, since otherwise $\Lambda (q_\lambda ,\nu ) \cap \lambda \neq \emptyset$. And $\lambda (f_{(0,q_\lambda ,\nu )}) \leq \lambda$, because $\lambda \in C_\nu$. Thus $q_\lambda \in rng(f)$ for some $f:\bar \nu \Rightarrow \nu$ by lemma 6 (b). But then $rng(f_{(0,q_\lambda ,\nu )} ) \subseteq rng(f)$. $\Box$
\medskip \\
{\bf Lemma 8}
\smallskip \\
Let $\rho \in C_\nu \cap \lambda$ such that $\rho > q_\lambda$. Then $q_\lambda$ is an initial segment of $q_\rho$.
\smallskip \\
{\bf Proof:}
\smallskip 

$q_\rho (k)=max(\Lambda (q_\rho \upharpoonright k ,\rho ))=max(\Lambda (q_\rho \upharpoonright k ,\nu ) \cap \rho )$,
\smallskip \\
as long as these maxima exist, because $\rho \in C_\nu$. Hence $q_\rho \upharpoonright k \in rng(f)$ for some $f:\bar \nu \Rightarrow \nu$ by lemma 6 (b). So $\Lambda (q_\rho \upharpoonright k ,\nu ) \cap \rho =\Lambda (q_\rho \upharpoonright k ,\rho )$ by lemma 5. Analogously
\smallskip 

$q_\lambda (k)=max(\Lambda (q_\lambda \upharpoonright k ,\lambda ))=max(\Lambda (q_\lambda \upharpoonright k ,\nu ) \cap \lambda )=max(\Lambda (q_\lambda \upharpoonright k ,\nu ) \cap \rho )$,
\smallskip \\
as long as these maxima exist, because $q_\lambda < \rho < \lambda$. The lemma follows from these two equations by induction. $\Box$
\bigskip \\
{\bf Lemma 9}
\smallskip \\
$C_\nu$ is closed in $\nu$.
\smallskip \\
{\bf Proof:} Let $\lambda \in Lim (C_\nu )$. Consider the sequence $\langle q_\rho \mid \rho \in C_\nu \cap \lambda \rangle$. By lemma 8, there is some $\rho _0 \in C_\nu \cap \lambda$ such that $q_\rho =q_{\rho _0}$ for all $\rho _0 < \rho \in C_\nu \cap \lambda$. Therefore, by lemma 7, $\rho \in \Lambda (q_{\rho _0},\nu )$ for all $\rho _0 < \rho \in C_\nu \cap \lambda$. But $\Lambda (q_{\rho _0},\nu )$ is closed. Hence $\lambda \in \Lambda (q_{\rho _0},\nu ) \subseteq C_\nu$. $\Box$
\medskip \\
{\bf Lemma 10}
\smallskip \\
$\lambda \in C_\nu \Rightarrow C_\lambda =C_\nu \cap \lambda$.
\smallskip \\
{\bf Proof} by induction on $\lambda$ and $\nu$. Suppose the lemma to be proved already for all $\rho < \lambda$ and $\mu \leq \nu$. By lemma 7, $\Lambda (q_\lambda ,\lambda )=\Lambda (q_\lambda ,\nu ) \cap \lambda$. Therefore $\rho \in C_\nu \cap C_\lambda$ for all $\rho \in \Lambda (q_\lambda , \lambda )$. Hence $C_\lambda \cap \rho =C_\nu \cap \rho =C_\rho$ by the induction hypothesis. If $\Lambda (q_\lambda ,\lambda )$ is unbounded in $\lambda$, we are finished. If $\Lambda (q_\lambda ,\lambda )=\emptyset$, then $(C_\nu \cap \lambda )-(q_\lambda (k_\lambda ) +1)=\emptyset$ by lemma 8. To see this, assume $(C_\nu \cap \lambda )-(q_\lambda (k_\lambda ) +1)\neq \emptyset$. Let $\rho =min(C_\nu -(q_\lambda (k_\lambda )+1))$. Then $q_\rho =q_\lambda$ by lemma 8. Hence $\rho \in \Lambda ( q_\lambda ,\lambda )$. Contradiction! Therefore $(C_\nu \cap \lambda )-(q_\lambda (k_\lambda )+1)=C_\lambda -(q_\lambda (k_\lambda )+1)=\emptyset$. If $q_\lambda (k_\lambda )=0$, then we are finished. If $q_\lambda (k_\lambda )\neq 0$, then $q_\lambda (k_\lambda )=max(C_\lambda )=max(C_\nu \cap \lambda )$. But $C_\lambda \cap q_\lambda (k_\lambda )=C_\nu \cap q_\lambda (k_\lambda )=C_{ q_\lambda (k_\lambda )}$. Hence $C_\lambda =C_\nu \cap \lambda$. $\Box$ 
\medskip \\
{\bf Lemma 11}
\smallskip \\
Let $f:\bar \nu \Rightarrow \nu$. Then $(f \upharpoonright J^D_{\bar \nu}) : \langle J^D_{\bar \nu} ,D \upharpoonright {\bar \nu},C_{\bar \nu}\rangle \rightarrow \langle J^D_\nu ,D\upharpoonright \nu ,C_\nu \rangle$ is $\Sigma _0$-elementary.
\smallskip \\
{\bf Proof:} Show $f(C_{\bar \nu} \cap \bar \eta )=C_\nu \cap f(\bar \eta )$ for all $\bar \eta < \bar \nu$. By (LP1), we have $f (C_{\bar \nu} \cap \bar \lambda ) = f (C_{\bar \lambda})=C_\lambda =C_\nu \cap f(\lambda )$ for all $\bar \lambda \in C_{\bar \nu}$. Therefore, if $C_{\bar \nu}$ is cofinal in $\bar \nu$, we are finished. If it is not, then $f(q_{\bar \nu})=q_\nu$. If $q_{\bar \nu}(k_{\bar \nu})=0$, then $\Lambda (0,\bar \nu )=\Lambda (0,\nu )=\emptyset$, implying that $C_{\bar \nu}=C_\nu =\emptyset$. If $q_{\bar \nu}(k_{\bar \nu}) \neq 0$, then we use $f(max(C_{\bar \nu}))=max(C_\nu )$. But $max(C_{\bar \nu})=q_{\bar \nu}(k_{\bar \nu})$ and $max(C_\nu )=q_\nu (k_\nu )$. $\Box$
\medskip \\
{\bf Lemma 12}
\smallskip \\
Set $\alpha _{\tau (0,\nu )}=\mu _\nu$ and $x (0,\nu )= \emptyset$ for all $\nu$. Then the following holds for all $0 \leq n$ and $\nu \in \widehat{S}$:
\smallskip \\
(i) If $f: \bar \nu \Rightarrow _{n+1} \nu$, $\alpha := \alpha _{\tau (n,\nu )}$ and $\bar \alpha := f^{-1}[\alpha \cap rng(f)]$, then $\bar \alpha =\alpha _{\tau (n,\bar \nu )}$.
\smallskip \\
(ii) If $f : \bar\nu \Rightarrow _{n+1} \nu$, then $f(x (n, \bar \nu ))=x (n ,\nu )$.
\smallskip \\
(iii) If $f:\bar \nu \Rightarrow _{n+1} \nu$ and $\bar K=f^{-1}[K^n_\nu \cap rng(f)]$, then $\bar K=K^n_{\bar \nu}$.
\smallskip \\
(iv) If $f,g\Rightarrow _{n+1} \nu$ and $rng(f) \subseteq rng(g)$, then $g^{-1}f \Rightarrow _{n+1} d(g)$.
\smallskip \\
(v) For all $u \subseteq J^D_{\mu _\nu}$, there is $f^{n+1}_{(u,\nu )}$.
\smallskip \\
(vi) For all $\beta < \nu$ and $x \in J^D_{\mu _\nu}$, $f^{n+1}_{(\beta , x ,\nu )}$ is uniformly definable over $\langle J^D_\nu ,D\upharpoonright \nu ,D_\nu \rangle$.
\smallskip \\
{\bf Proof} by induction on $n$. For $n=0$, (i) to (v) hold by the morass axioms. 
\smallskip \\
(vi) By (DF), the $rng(f^1_{(0,x ,\nu )})$ are uniformly definable over $\langle J^D_\nu ,D\upharpoonright\nu ,D_\nu \rangle$. Like in the proof of lemma 1, $rng(f^1_{(\beta ,z_0 ,\nu )})=\bigcup \{ rng(f^1_{(0,z ,\nu )}) \mid z \in (\beta \cup \{ z_0 \} )^{<\omega} \}$. And $f^1_{( \beta ,z_0 ,\nu )}(x)=y$ may be defined by: There is some $\bar \nu$ and some $\bar z_0$ such that, for all $z_1 \in \beta ^{<\omega}$,
$$d(f_{(0,\langle z_1 ,\bar z_0 \rangle ,\bar \nu )})=d(f_{(0,\langle z_1 ,z_0 \rangle ,\nu )})$$
$and$, for all $t \in J^D_{\bar \nu}$, there is some $z_1 \in \beta ^{<\omega}$ such that
$$t \in rng(f_{(0,\langle z_1 ,\bar z_0 \rangle ,\bar \nu )})$$
$and$ there is some $z$ and some $z_1 \in \beta ^{<\omega}$ such that
$$f_{(0,\langle z_1 ,\bar z_0 \rangle ,\bar \nu )}(z)=x \Leftrightarrow f_{(0,\langle z_1 ,z_0 \rangle ,\nu )}(z)=y .$$  
Now, assume that (i) to (vi) are proved already for all $0 \leq m <n$.
\smallskip \\
(i) Let $B^n(x ,\nu ):=\{ \beta (f^n_{(\gamma , x ,\nu )})<\alpha _{\tau (n,\nu )} \mid \gamma < \nu \} = \{ \beta < \alpha _{\tau (n,\nu )}\mid \beta \notin rng(f^n_{(\beta ,x ,\nu )})\}$. Let $f(\bar x )=x (n,\nu )$, $B=B^n(x ,\nu )$ and $\bar B:= f^{-1}[B \cap rng(f)]$. Then $f \circ f^n_{(\bar u,\bar \nu )}=f^n_{(u,\nu )}$ for all $\bar u \subseteq J^D_{\mu _{\bar\nu}}$ and $u=f[\bar u]$ by (iv) of the induction hypothesis (cf. lemma 2b). Therefore, if $f(\bar \beta )=\beta \in rng(f_{(\beta ,x (n,\nu ),\nu )})$, then $\bar \beta \in rng(f_{(\bar \beta ,\bar x ,\bar \nu )})$, because $\bar \beta = f^{-1}[\beta \cap rng(f)]$. And if $f(\bar \beta )=\beta \not\in rng(f_{(\beta ,x (n,\nu ),\nu )})$, then $\bar \beta \not\in rng(f_{(\bar \beta ,\bar x ,\bar \nu )})$. So, altogether, $\bar B=B^n(\bar x ,\bar \nu )$. By (DP3)(b) and (iv) of the induction hypothesis, $B^n(x (n,\nu ),\nu )=\bigcup \{ B^n(x (n,\eta ),\eta )\mid \eta \in K^n_\nu \}$. But $B^n(x (n,\nu ),\nu )$ is unbounded in $\alpha$ and $rng(f) \cap J^D_\alpha \prec _1 \langle J^D_\alpha ,D\upharpoonright \alpha ,K^n_\nu \rangle$. Thus $\bar B=B^n(\bar x,\bar\nu )$ is also unbounded in $\bar\alpha$. Assume there were some $z \in J^D_{\mu _{\bar\nu}}$ and some $\beta < \bar \alpha$ such that $f^n_{(\beta , z ,\bar \nu )}=id _{\bar \nu}$. Then there was some $\beta \leq \gamma < \bar \alpha$ such that $z \in rng(f^n_{(\gamma , \bar x ,\bar \nu )})$. For, by (iv) of the induction hypothesis, $f^n_{(\bar \alpha ,\bar x ,\bar \nu )}=id_{\bar\nu}$. So $f^n_{(\gamma , \bar x ,\bar \nu )}=id _{\bar \nu}$. But this contradicts the fact that $B^n(\bar x ,\bar \nu )$ is unbounded in $\bar \alpha$.  
\smallskip \\
(ii) By the proof of (i), $f^n_{(\bar \alpha ,\bar x ,\bar \nu )}=id_{\bar \nu}$ is satisfied for  
$\bar \alpha =\alpha _{\tau (n,\bar \nu )}$ and $f(\bar x )=x (n,\bar \nu )$. Therefore $x (n,\bar \nu ) \leq \bar x$. Assume $x (n,\bar \nu ) < \bar x$. Then $x (n,\nu )\in rng(f^n_{(\alpha ,x ,\nu )})$ where $x := f(x (n,\bar \nu ))$ and $\alpha := \alpha _{\tau (n,\nu )}$. Thus $f^n_{(\alpha ,x ,\nu )}=id_\nu$ for all $x < x (n,\nu )$. But that contradicts the definition of $x (n,\nu )$.
\smallskip \\
(iii) Let $f(\bar \mu )=\mu$, $K^+=K^n_\nu -Lim (K^n_\nu )$ and $\bar K^+=K^n_{\bar \nu} -Lim (K^n_{\bar \nu} )$. First, prove $\mu \in K^+ \Rightarrow \bar \mu \in \bar K^+$. By (i) and (ii), we know that $\bar B=f^{-1}[B \cap rng(f)]$ where $B=B^n(x ,\nu )$, $\bar B=B^n(\bar x ,\bar \nu )$ and $x =x (n,\nu )$, $\bar x =x (n,\bar \nu )$. Let $\mu =d(f^n_{(\beta ,x ,\nu )})$. Since $\mu \in K^+ \cap rng(f)$, we may assume $\beta \in B^+ \cap rng(f)$ where $B^+=B-Lim(B)$. Let $\delta$ be the predecessor of $\beta$ in $B$. Then $f^n_{(\delta , \langle \delta ,x \rangle ,\mu )}=id_\mu$. Define $\gamma = \beta$ if $\beta \in S^+ \cup S^0$, and $\gamma = min \{ \gamma \sqsubset \beta \mid \delta < \gamma \}$ else. Then $\gamma \in rng(f)$ and $\mu =\mu _\gamma$ by (DP3). Let $f(\bar \beta )=\beta$, $f(\bar \gamma )=\gamma$. By (iv) of the induction hypothesis, $\bar \mu = \mu _{\bar\gamma}=d(f^n_{(\bar \beta ,\bar x ,\bar \nu )}) \in \bar K^+$. In the same way, we show $\bar \mu \in \bar K^+ \Rightarrow \mu \in K^+$. But $K^n_\nu = \bigcup \{ K^n_\eta \mid \eta \in K^+ \}$ and $K^n_{\bar\nu} = \bigcup \{ K^n_\eta \mid \eta \in \bar K^+ \}$. Thus the claim holds.
\smallskip \\
(iv) follows immediately from (ii), (iii) and the definition of $\Rightarrow _{n+1}$.
\smallskip \\
(v) First, we notice that $\langle J^D_\alpha ,D\upharpoonright \alpha ,K^n_\nu \rangle$ where $\alpha := \alpha _{\tau (n,\nu )}$ is rudimentary closed. Then $K^n_\nu \cap \eta = K^n_\eta$ for all $\eta \in K^n_\eta$ by (iv). But, by (vi) of the induction hypothesis, $K^n_\eta$ is uniformly definable over $\langle J^D_\eta ,D\upharpoonright \eta ,D_\eta \rangle$. Since $\langle J^D_\alpha ,D\upharpoonright \alpha ,K^n_\nu \rangle$ is  rudimentary closed, by the definition of $\Rightarrow _{n+1}$,
$$f^{n+1}_{(u,\nu )}=f^n_{(w \cup u \cup \{ x (n,\nu )\} ,\nu )}$$

where $w:= h[\omega \times (u \cap J^D_\alpha )^{<\omega}]$.
\smallskip \\
Here, $h$ denotes the canonical $\Sigma _1$-Skolem function of $\langle J^D_\alpha ,D\upharpoonright\alpha ,K^n_\nu \rangle$. 
\smallskip \\
(vi) If $w \prec _1 \langle J^D_{\alpha _{\tau (n,\nu )}},D\upharpoonright {\alpha _{\tau (n,\nu )}},K^n_\nu \rangle$, then there is a uniquely determined $f \Rightarrow _{n+1} \nu$ such that $rng(f)\cap J^D_{\alpha _{\tau (n,\nu )}}=w$.
\smallskip   \\
Existence :
\smallskip \\
Let $\alpha := \alpha _{\tau (n,\nu )}$ and
$$f_\beta =f^n_{(\beta , x (n,\nu ),\nu )}$$
$$\nu (\beta )=d(f_\beta )$$
$$H= \bigcup \{ f_\beta [w \cap J^D_{\nu (\beta )}]\mid \beta < \alpha \} .$$  
Then $H \cap J^D_\alpha =w$. For $w \subseteq H \cap J^D_\alpha$ is clear, since $f_\beta \upharpoonright J^D_\beta = id \upharpoonright J^D_\beta$. So let $y \in H \cap J^D_\alpha$. Thus $y=f_\beta (x)$ for some $x \in w$ and some $\beta < \alpha$. Let 
$K^+ =K^n_\nu -Lim(K^n_\nu )$ and $\beta (\eta )= sup \{ \beta \mid f^n_{(\beta ,x (n,\eta ), \eta )}\neq id_\eta \}$. Then
$$\langle J^D_\alpha , D\upharpoonright \alpha ,K^n_\nu \rangle \models (\exists y)(\exists \eta \in K^+ )(y=f^{m+1}_{(\beta ,x (m+1,\eta ), \eta )}(x) \in J^D_{\beta (\eta )}).$$
Since $w \prec _1 \langle J^D_\alpha ,D\upharpoonright\alpha ,K^n_\nu \rangle$, $y=f^n_{(\beta ,x (n,\eta ), \eta )}(x)$ $\in w$ for all such $\eta$ and $x \in w$. But since $y=f^n_{(\beta ,x (n,\eta ), \eta )}(x) \in J^D_{\beta (\eta )}$, we get $f_\beta (x) = f^n_{(\beta ,x (n,\eta ), \eta )}(x) \in w$.

Let $|f| : J^D_{\bar \nu} \rightarrow J^D_\nu$ be the uncollapse of $H$ and $f = \langle \bar \nu ,|f | ,\nu \rangle$. Then $f:\bar \nu \Rightarrow _{n+1} \nu$. For, for all $\beta < \alpha$ by (DF), $f^{(\bar \nu (\beta ))}:\bar \nu (\beta ) \Rightarrow _n \nu (\beta )$ where $f(\bar \nu (\beta ))=\nu (\beta )$ if $\nu (\beta ) \in rng(f)$. Let $\Gamma = \{ \beta < \alpha  \mid \nu (\beta ) \in rng(f) \}$. For $\beta ,\gamma \in \Gamma$, let $g_\beta =f_\beta \circ f^{(\bar \nu (\beta ))}$ and $g_{\beta \gamma}=g^{-1}_\gamma \circ g_\beta$. Let $\langle h_\beta \mid \beta \in \Gamma \rangle$ be the transitive, direct limit of the directed system $\langle g_{\beta \gamma}\mid \beta \leq \gamma \in \Gamma \rangle$. Then $f \circ h_\beta = g_\beta$ for all $\beta \in \Gamma$. Thus, by (CP1) and (iv) of the induction hypothesis, $f:\bar \nu \Rightarrow _n \nu$. But $x (n+1,\nu )\in H=rng(f)$ and $rng(f) \cap J^D_\alpha = w \prec _1 \langle J^D_\alpha ,D\upharpoonright\alpha ,K^n_\nu \rangle$. Thus $f:\bar \nu \Rightarrow _{n+1} \nu$.  
\smallskip \\
Uniqueness:
\smallskip \\
Let $f:\bar \nu \Rightarrow _{n+1} \nu$ such that $rng(f)\cap J^D_{\alpha _{\tau (n,\nu )}}=w$ and $\bar \alpha := f^{-1}[\alpha \cap rng(f)]$. Then $\bar \alpha = \alpha _{\tau (n,\bar\nu )}$ by (i). And $f \circ f^{n+1}_{(\bar \alpha , \bar \nu )}=f^{n+1}_{(w,\nu )}$ by (iv) (cf. lemma 2a). But $f^{n+1}_{(\bar \alpha ,\bar \nu )}=id_{\bar\nu}$, since $\bar \alpha = \alpha _{\tau (n,\bar \nu )}$. Therefore, $f=f^{n+1}_{(w,\nu )}$ is uniquely determined. 
\smallskip \\
Let $f^n_{(0,\langle x (n,\nu ),z_0 \rangle ,\nu )}(z_0 ^\ast )=z_0$. Use $w=h_{(n,\nu )}[\omega \times (\beta ^{<\omega} \times \{ z_0 \} )]$ where $h_{(n,\nu )}$ is the canonical $\Sigma _1$-Skolem function of $\langle J^D_{\alpha _{\tau (n,\nu )}} ,D\upharpoonright{\alpha _{\tau (n,\nu )}},K^n_\nu \rangle$. By (vi) of the induction hypothesis, $K^n_\nu$ is  uniformly definable over $\langle J^D_\nu ,D\upharpoonright\nu ,D_\nu \rangle$. Therefore, $w$ is uniformly definable over $\langle J^D_\nu ,D\upharpoonright\nu ,D_\nu \rangle$. Let $\pi$ be the uncollapse of $w$. Then we can define $\pi (x)=y$ by: 
There is some $\bar \nu \leq \nu$ and some $\bar z_0 \leq z_0 ^\ast$ such that, for all $i \in \omega$ and $z_1 \in \beta ^{<\omega}$,
$$(\exists z \in J^D_{\alpha _{\tau (n,\bar\nu )}})(z=h_{(n,\bar\nu )}(i,\langle z_1 ,\bar z_0 \rangle )) \Leftrightarrow  (\exists z \in J^X_{\alpha _{\tau (n,\nu )}})(z=h_{(n,\nu )}(i,\langle z_1 ,z_0 ^\ast \rangle ))$$
$and$, for all $z \in J^X_{\alpha _{\tau (n,\bar\nu )}}$, there is some $i \in \omega$ and some $z_1 \in \beta ^{<\omega}$ such that
$$z=h_{(n,\bar\nu )}(i,\langle z_1 ,\bar z_0 \rangle )$$
$and$ there is some $i \in \omega$ and some $z_1 \in \beta ^{<\omega}$ such that
$$h_{(n,\bar \nu )}(i,\langle z_1 ,\bar z_0 \rangle )=x \Leftrightarrow h_{(n,\nu )}(i,\langle z_1 ,z_0 ^\ast \rangle )=y .$$ By this, $\bar \nu$ is uniquely determined. By what was shown above, one can define $f^{n+1}_{(\beta ,z_0 ,\nu )}(x)=f^n_{(w,\nu )}(x)=y$ by: For all $z_0 \in \alpha _{\tau (n,\bar \nu )} ^{<\omega}$,
$$d(f^n_{(0,\langle z_0 ,x (n,\bar \nu ) \rangle ,\bar \nu )})=d(f^n_{(0,\langle \pi (z_0 ) ,x (n,\nu ) \rangle ,\nu )})$$
$and$, for all $t \in J^D_{\bar \nu}$, there is some $z_0 \in \alpha _{\tau (n,\bar \nu )} ^{<\omega}$ such that
$$t \in rng(f^n_{(0,\langle z_0 ,x (n, \bar \nu ) \rangle ,\bar \nu )})$$
$and$ there is some $z$ and some $z_0 \in \alpha _{\tau (n,\bar \nu )} ^{<\omega}$ such that
$$f^n_{(0,\langle z_0 ,x (n,\bar \nu ) \rangle ,\bar \nu )}(z)=x \Leftrightarrow f^n_{(0,\langle \pi (z_0 ) ,x (n,\nu ) \rangle ,\nu )}(\pi (z))=y .$$  
$\Box$
\medskip\\
Now, it is an immediate consequence of lemma 12 and (DF) that ($\times$) holds for all $\nu$ such that $\mu _\nu =\mu$.

\section{The inner model $L[X]$}
Of course my definition of $(\omega _1,\beta)$-morass makes also sense if $\beta < \omega_1$. Hence a natural question is:
\smallskip\\
Is the existence of an $(\omega _1, \beta)$-morass in this new sense equivalent to the existence of an $(\omega _1, \beta)$-morass in Jensen's sense?
\smallskip\\
In asking this question one has to be careful what an $(\omega _1, \beta)$-morass in Jensen's sense is, because there are also different definitions. But for the case $\beta =1$, I expect an equivalence between all existing definitions.
\smallskip\\
In the following, I will define a strengthening of the notion of a Jensen $(\omega _1, \beta)$-morass which I also expect to be equivalent to my notion of $(\omega _1,\beta)$-morass. If we construct a morass in the usual way in $L$, the properties of this stronger notion hold automatically (see the paper [Irr2] or my dissertation [Irr1]).
\medskip\\
A structure  $\frak{M} = \langle S,\vartriangleleft ,\frak{F},D \rangle$ is called an $\omega _{1+\beta}$-standard morass if it satisfies all axioms of an $( \omega _1 , \beta )$-morass except {\bf (DF)} which is replaced by:
\medskip 

$\nu \vartriangleleft \tau$ $\Rightarrow$ $\nu$ is regular in $J^D_\tau$
\medskip \\
and there are functions $\sigma _{(x ,\nu )}$ for $\nu \in \widehat{S}$ and $x \in J^D_\nu$ such that:
\pagebreak \\
{\bf (MP)$^+$}
\smallskip \\
$\sigma _{(x ,\nu )}[\omega ]=rng(f_{(0,x ,\nu )})$
\medskip \\
{\bf (CP1)$^+$}
\smallskip \\
If $f:\bar \nu \Rightarrow \nu$ and $f(\bar x )=x$, then $\sigma _{(x ,\nu )}=f \circ \sigma _{(\bar x ,\bar \nu )}$.
\medskip \\
{\bf (CP3)$^+$}
\smallskip \\
If $C_\nu$ is unbounded in $\nu$, then $\sigma _{(x ,\nu )}=\bigcup \{ \sigma _{(x ,\lambda )} \mid \lambda \in C_\nu , x \in J^D_\lambda \}$.
\medskip \\
{\bf (DF)$^+$}
\smallskip \\
(a) If $f_{(0,x ,\nu )}=id_\nu$ for some $x \in J^D_\nu$, then
$$\{ \langle i ,z ,\sigma _{(z ,\nu )}(i) \rangle \mid z \in J^D_\nu, i \in dom(\sigma _{(z ,\nu )}) \}$$
is uniformly definable over $\langle J^D_{\mu _\nu},D\upharpoonright {\mu _\nu},D_{\mu _\nu} \rangle$.
\smallskip \\
(b) If $C_\nu$ is unbounded in $\nu$, then $D_\nu =C_\nu$. If it is bounded, then $D_\nu =\{ \langle i,\sigma _{(q_\nu ,\nu )} (i) \rangle \mid i \in dom (\sigma _{(q_\nu ,\nu )}) \}$. 
\medskip \\
{\bf Lemma 13}
\smallskip \\
(DF) and ($\times$) also hold in a standard morass.
\smallskip \\
{\bf Proof:} First, we prove by induction on $\mu \in \widehat{S}$ that the set  
$$\{ \langle i ,x ,\sigma _{(x ,\mu )}(i) \rangle \mid x \in J^D_\mu, i \in dom(\sigma _{(x ,\mu )}) \}$$
is uniformly definable over $\langle J^D_\mu ,D\upharpoonright\mu ,D_\mu \rangle$ for all  $\mu \in \widehat{S}$.
Assume that this has been proved already for all $\tau <\mu$, $\tau \in \widehat{S}$. 
\smallskip \\
If $C_\mu$ is unbounded in $\mu$, then, by (CP3)$^+$, 
$$\sigma _{(x ,\nu )}=\bigcup \{ \sigma _{(x ,\lambda )} \mid \lambda \in C_\nu , x \in J^D_\lambda \} .$$ 
But, by the induction hypothesis, the $\sigma _{(x ,\lambda )}$, $\lambda \in C_\mu$, are uniformly definable over $\langle J^D_\mu ,D\upharpoonright\mu ,D_\mu \rangle$. And, by (DF)$^+$(b), $C_\mu = D_\mu$. Thus $\sigma _{(x ,\nu )}$ is  uniformly definable over $\langle J^D_\mu ,D\upharpoonright\mu ,D_\mu \rangle$. 
\smallskip \\
If $C_\mu$ is bounded in $\mu$, then $rng(f_{(0,q_\mu ,\mu )})$ is unbounded in $\mu$. Therefore, by (CP2),
$$rng(f_{(0,\langle z_0 ,q_\mu \rangle ,\mu )})=h_\mu [\omega \times (rng(f_{(0,q_\mu ,\mu )}) \times \{ z_0 \} )].$$
Here, $h_\mu$ is the $\Sigma _1$-Skolem function of $\langle J^D_\mu ,D\upharpoonright \mu \rangle$. Since $D_\mu =rng(f_{(0,q_\mu ,\mu )})$, the $rng(f_{(0,\langle z_0 ,q_\mu \rangle ,\mu )})$ are uniformly definable over $\langle J^D_\mu ,D\upharpoonright\mu ,D_\mu \rangle$. Since, by (CP1)$^+$ and lemma 6 (b), for $\bar \mu := d(f_{(0,\langle z_0 ,q_\mu \rangle ,\mu )})$, 
$$f_{(0,\langle z_0 ,q_\mu \rangle ,\mu )}\circ \sigma _{(q_{\bar\mu} ,\bar\mu )} = \sigma _{(q_\mu ,\mu )}$$
holds, we can define $f_{(0,\langle z_0 ,q_\mu \rangle ,\mu )}(x)=y$ by: 
There is some $\bar \mu \leq \mu$ and some $\bar z_0 \leq z_0$ such that, for all $i,j \in \omega$,
$$(\exists z \in J^D_{\bar \mu})(z=h_{\bar \mu }(i,\langle \sigma _{(q_{\bar\mu} ,\bar \mu )}(j) , \bar z_0 \rangle )) \Leftrightarrow  (\exists z \in J^D_{\mu})(z=h_\mu (i,\langle \sigma _{(q_{\mu} ,\mu )}(j),z_0 \rangle ))$$
$and$, for all $z \in J^D_{\bar\mu}$, there is some $i \in \omega$ and some $j \in \omega$ such that
$$z=h_{\bar\mu}(i,\langle \sigma _{(q_{\bar\mu} ,\bar\mu )}(j),\bar z_0 \rangle )$$
 $and$ there is some $i \in \omega$ and some $j\in \omega$ such that
$$h_{\bar\mu}(i,\langle \sigma _{(q_{\bar\mu} ,\bar\mu )}(j),\bar z_0 \rangle )=x \Leftrightarrow h_\mu (i,\langle \sigma _{(q_{\mu} ,\mu )}(j),z_0 \rangle )=y .$$
\smallskip \\
If $\alpha _{\tau (1,\mu )}=0$, then it follows from (DF)$^+$ that $\{ \langle i ,z_0 ,\sigma _{(z_0 ,\mu )}(i) \rangle \mid z_0 \in J^D_\mu, i \in dom(\sigma _{(z_0 ,\mu )}) \}$ is uniformly definable over $\langle J^D_\mu ,D\upharpoonright\mu ,D_\mu  \rangle$. If $\alpha _{\tau (1,\mu )} >0$, then, by (DP3)(b), $\bar \mu = d(f_{(0,\langle z_0 ,q_\mu \rangle ,\mu )}) < \mu$. But then, by (CP1)$^+$, $\sigma _{(z_0 ,\mu )}=f_{(0,\langle z_0 ,q_\mu \rangle ,\mu )} \circ \sigma _{(\bar z_0 ,\bar \mu )}$, where $f_{(0,\langle z_0 ,q_\mu \rangle ,\mu )}(\bar z_0 )=z_0$, is definable by the induction hypothesis.
\smallskip \\ 
From the $\sigma$s, we calculate $f_{(0,z_0 ,\mu )}(x)=y$ as follows:
There is some $\bar \mu \leq \mu$ and some $\bar z_0 \leq z_0$ such that, for all $r,s \in \omega$,
$$\sigma _{(\bar z_0 ,\bar \mu )}(r) \leq \sigma _{(\bar z_0 ,\bar \mu )}(s) \Leftrightarrow \sigma _{(z_0 ,\mu )}(r) \leq \sigma _{(z_0 ,\mu )}(s)$$
$and$, for all $z \in J^D_{\bar \mu}$, there exists some $s \in \omega$ such that
$$z=\sigma _{(\bar z_0 ,\bar \mu )}(s)$$
 $and$ there exists some $s \in \omega$ such that
$$\sigma _{(\bar z_0 ,\bar \mu )}(s)=x \Leftrightarrow \sigma _{(z_0 ,\mu )}(s)=y .$$
Since the $f^1_{(0,z_0 ,\mu )}$ are uniformly definable over $\langle J^D_\mu ,D\upharpoonright\mu ,D_\mu  \rangle$ and (DF) and ($\times$) hold by the induction hypothesis for all $\tau \in \widehat{S} \cap \mu$, we can define the $f^n_{(0,z_0 ,\mu )}$ with $z_0 \in J^D_\mu$ uniformly over $\langle J^D_\mu ,D\upharpoonright\mu ,D_\mu  \rangle$ like in the proof of lemma 12. Finally,  $$\{ \langle z_0 ,\nu ,x,f_{(0,z_0 ,\nu )}(x) \rangle \mid  \nu < \mu ,\mu _\nu =\mu , z_0 \in J^D_\mu , x \in dom(f_{(0,z_0 ,\nu )}) \}$$
$$\cup \{ \langle z_0 ,x,f_{(0,z_0 ,\mu )}(x) \rangle \mid z_0 \in J^D_\mu , x \in dom(f_{(0,z_0 ,\mu )}) \}$$
$$\cup (\sqsubset \cap \mu ^2)$$
may be defined over $\langle J^D_\mu ,D\upharpoonright\mu ,D_\mu \rangle$ using (DF).
$\Box$
\medskip\\
Let $S^X \subseteq Lim$ and $X=\langle X_\nu \mid \nu \in S^X \rangle$ be a sequence.
\smallskip \\
Let $I_\nu =\langle J^X_\nu , X \upharpoonright \nu \rangle$ for $\nu \in Lim -S^X$ and $I_\nu = \langle J^X_\nu , X \upharpoonright \nu ,X_\nu \rangle$ for $\nu \in S^X$ where $X_\nu \subseteq J^X_\nu$ and
\smallskip

$J^X_0 = \emptyset$
\smallskip

$J^X_{\nu +\omega } = rud(I^X_\nu )$
\smallskip

$J^X_\lambda = \bigcup \{ J^X_\nu \mid \nu \in \lambda \}$ for $\lambda \in Lim^2:=Lim(Lim)$.
\smallskip \\
Here, $rud(I^X_\nu )$ is the rudimentary closure of $J^X_\nu \cup \{ J^X_\nu \}$ relative to $X \upharpoonright \nu$ if $\nu \in Lim - S^X$ and relative to $X \upharpoonright \nu$ and $X_\nu$ if $\nu \in S^X$.
\medskip\\
Let $\beta (\nu )$ be the least $\beta$ such that $J^X_{\beta +\omega} \models \nu$ singular.
\medskip\\
Now, let a $\kappa$-standard morass be given. I will show that there is an $S^X\subseteq \kappa$ and a sequence $X$ as above that the following holds:
\smallskip\\
{\bf (Amenability)} The structures $I_\nu$ are amenable.
\smallskip\\
{\bf (Coherence)} If $\nu \in S^X$, $H \prec _1 I_\nu$ and $\lambda = sup(H \cap On)$, then $\lambda \in S^X$ and $X_\lambda = X_\nu \cap J^X_\lambda$.
\smallskip \\
{\bf (Condensation)} If $\nu \in S^X$ and $H \prec_1 I_\nu$, then there is some $\mu \in S^X$ such that $H \cong I_\mu$.
\smallskip\\
{\bf ($\ast$)} $Card \cap \kappa = Card^{L_\kappa [X]}$.
\smallskip\\
{\bf ($\ast\ast$)} $S^X=\{ \beta (\nu )\mid \nu$ singular in $I_\kappa \}$.
\medskip \\
These properties are good enough to do fine structure proofs in $L_\kappa [X]$, e.g. to construct a $\kappa$-standard morass. This will be shown in a forthcoming paper [Irr2].
\medskip\\
To define $X$, I will use the sets $C_\nu$ from (CP3):
\medskip \\
If $\nu \in \widehat{S}$ and $C_\nu$ is unbounded in $\nu$, then set 
$$X_\nu =C_\nu .$$
Let $\nu \in \widehat{S}$ and $C_\nu$ be bounded in $\nu$. Then $\Lambda (q,\nu )$ is bounded for all $q \in \nu$. Thus $\Lambda (q_\nu ,\nu )=\emptyset$. So $f_{(0,q_\nu ,\nu )}$ is cofinal. In this case, set
$$X_\nu =\{\sigma _{(q_\nu ,\nu )}  [n] \mid n \in \omega \} .$$
\smallskip \\
Let $S^X =\widehat{S}$.
\medskip \\
{\bf Lemma 14}
\smallskip \\
If $\nu \in \widehat{S}$, $C_\nu$ is unbounded in $\nu$ and
$f:\langle J^{\bar D}_{\bar \nu} ,\bar D,\bar C \rangle \rightarrow \langle J^D_\nu ,D\upharpoonright\nu ,C_\nu \rangle$ is $\Sigma _1$-elementary, then $\langle \bar \nu , f ,\nu \rangle \in \frak{F}$.
\smallskip \\
{\bf Proof:} Let $z_0 \in rng(f)$, $i \in \omega$ and $y=\sigma _{(z_0 ,\nu )}(i)$. Then we must prove $y \in rng(f)$. Since $C_\nu$ is unbounded in $\nu$, there is some $\lambda \in C_\nu$ such that $y=\sigma _{(z_0 ,\lambda )}(i)$ by (CP3)$^+$. Since, by lemma 13, the $\sigma _{(z_0 ,\tau )}$ are definable in $\langle J^D_\nu ,D\upharpoonright\nu \rangle$ when $\tau < \nu$, we have $\langle J^D_\nu , D\upharpoonright\nu ,C_\nu \rangle \models (\exists y)(\exists \lambda \in C_\nu )(y=\sigma _{(z_0 ,\lambda )}(i))$. Therefore, also $rng(f)\models (\exists y)(\exists \lambda \in C_\nu )(y=\sigma _{(z_0 ,\lambda )}(i))$. Thus $y \in rng(f)$. $\Box$
\medskip \\
{\bf Lemma 15}
\smallskip \\
Let $\nu \in \widehat{S}$, $H \prec _1 I_\nu$ and $f$ be the uncollapse of $H$. Let $f \upharpoonright On :\bar \nu \rightarrow \nu$. Then $\langle \bar \nu , f,\nu \rangle \in \frak{F}$.
\smallskip \\
{\bf Proof:} If $C_\nu$ is bounded in $\nu$, then $\lambda (f_{(0,q_\nu ,\nu )})=\nu$ and $rng(f_{(0,q_\nu ,\nu )}) \subseteq rng(f)$ by the definition of $X_\nu$. In addition, $f \upharpoonright J^D_{\bar \nu} : \langle J^D_{\bar \nu} , D\upharpoonright {\bar \nu}\rangle \rightarrow \langle J^D_\nu ,D\upharpoonright\nu \rangle$ is $\Sigma _1$-elementary. So the claim follows from (CP2). If $C_\nu$ is unbounded in $\nu$, then it follows from lemma 14. $\Box$
\medskip\\
{\bf Lemma 16}
\smallskip \\
(Coherence), (Amenability), (Condensation), ($\ast$) and ($\ast\ast$) hold for the sequence $X=\langle X_\nu \mid \nu \in S^X \rangle$. 
\smallskip \\
{\bf Proof:}
\smallskip \\
{\bf (Coherence)}
\smallskip \\
Let $\nu \in S^X$ and $H \prec _1 I_\nu$. If $C_\nu$ is unbounded in $\nu$, then $\lambda := sup(H \cap \nu ) \in C_\nu$ and $C_\nu \cap \lambda$ is unbounded in $\lambda$ by lemma 15. But, by lemma 10, $C_\nu \cap \lambda =C_\lambda$. So $X_\lambda = X_\nu \cap \lambda$. But if $C_\nu$ is bounded in $\nu$, then $H \cap \nu$ is unbounded in $\nu$ by the definition of $X_\nu$. So there is nothing to prove. 
\pagebreak \\
{\bf (Amenability)}
\smallskip \\
If $\nu \in S^X$ and $C_\nu$ is bounded in $\nu$, then $X_\nu \cap J^X_\eta$, $\eta < \nu$, is always finite. Therefore amenability is trivial. If $\nu \in S^X$ and $C_\nu$ is unbounded in $\nu$, then $C_\nu \cap \lambda =C_\lambda$ for all $\lambda \in C_\nu$ by lemma 10. Therefore $X_\lambda =X_\nu \cap \lambda$ for all $\lambda \in Lim(C_\nu )$ . If $Lim(C_\nu )$ is unbounded in $\nu$, we are finished. If it is not, then let $\lambda := max(Lim(C_\nu ))$. Then $X_\nu \cap J^X_\eta =C_\lambda \cup E$ where $E$ is finite for all $\eta > \lambda$.
\smallskip \\
{\bf(Condensation)}
\smallskip \\
If $\nu \in S^X$, $H \prec _1 I_\nu$ and $C_\nu$ is unbounded in $\nu$, then condensation holds by lemmas 11 und 15.
If $\nu \in S^X$ and $C_\nu$ is bounded in $\nu$, then $H \prec _1 I_\nu$ is unbounded in $\nu$ by the definition of $X_\nu$. Let $\pi$ be the uncollapse of $H$ and $\pi \upharpoonright On :\bar \nu \rightarrow \nu$. By lemmas 6 (b) and 15, $\pi (q_{\bar \nu})=q_\nu$. By the properties of $\sigma _\nu$ and $\sigma _{\bar \nu}$, we have condensation.
\smallskip \\
{\bf ($\ast$)}
\smallskip \\
Let $\omega < \kappa$ be a cardinal. Then all $\nu \in S_\kappa$ are independent by (DP1). Therefore $\bar \nu < \alpha _\nu =\kappa$ for all $f_{(\beta ,0,\nu )}:\bar \nu \Rightarrow \nu$ where $\beta < \alpha _\nu =\kappa$. Thus $rng(F)=\bigcup \{ rng(f_{(\beta ,0,\nu )}) \cap \nu \mid \beta < \alpha _\nu \} = \nu$ for $F:\{ \langle \beta ,x \rangle \mid x < d(f_{(\beta ,0,\nu )})\} \rightarrow \nu$ where $F(\beta ,x)=f_{(\beta ,0,\nu )}(x)$. By lemma 13, $F \in L_\kappa [X]$. So there is a map from a subset of $\kappa \times \kappa$ onto $\nu$ in $L_\kappa [X]$. By axioms (c) and (e), $S_\kappa$ is unbounded in $\kappa ^+$. Thus $(\kappa ^+)^{L_\kappa [X]}=\kappa ^+$. Since $\omega < \kappa$ was arbitrary, we get $Card ^{L[X]} - \omega _1=Card -\omega _1$. It remains to prove $\omega _1^{L_\kappa [X]}=\omega _1$. Let $\nu \in S_{\omega _1}$ and $\eta < \omega _1$. By axiom (1), $\eta \subseteq rng(f_{(0,\eta ,\nu )})$. By the definition of $X$, there exists a map from $\omega$ onto $\eta \subseteq rng(f_{(0,\eta ,\nu )})$ in $L_\kappa [X]$. If $n_\nu =1$, then $\sigma _{(\langle \eta , \alpha _\nu ^\ast ,P_\nu \rangle ,\mu _\nu )}$ is a map as needed by (DF). If $n_\nu > 1$,
$$h(i) := h_{\alpha _{\tau (n_\nu -1, \mu _\nu )},K^{n_\nu -1}_{\mu _\nu}}(i, \langle \eta ,\nu ^\ast ,\alpha _\nu ^{\ast\ast},P^\ast _\nu \rangle )$$
is as needed, by lemma 12 (vi) and (DF), where
\smallskip 

$f^{n_\nu -1}_{(\beta ,\langle x (n_\nu -1,\mu _\nu ), \alpha ^\ast _\nu \rangle ,\mu _\nu )}(\alpha ^{\ast\ast}_\nu )=\alpha ^\ast _\nu$
\smallskip 

$f^{n_\nu -1}_{(\beta , \langle x (n_\nu -1,\mu _\nu ) , P_\nu\rangle ,\mu _\nu )}(P^\ast _\nu )=P_\nu$
\smallskip 

$\nu ^\ast = \nu$ if $\nu < \alpha _{\tau (n_\nu -1, \mu _\nu )}$ and $\nu ^\ast = 0$ else.
\smallskip \\
Since $\eta < \omega _1$ was arbitrary, $\omega _1^{L_\kappa [X]}=\omega _1$. 
\smallskip \\
{\bf ($\ast \ast$)}
\smallskip \\
On the one hand, by definition of $n_\nu$ in (DF), there exists some $z_0 \in J^D_{\mu _\nu}$ and some $\gamma \sqsubset \nu$ such that $f^{n_\nu}_{(\gamma ,z_0 ,\mu _\nu )}$ is cofinal in $\nu$. If $n_\nu =1$, then $F:\gamma \times \omega \rightarrow \mu _\nu$ where
$$\langle \eta ,i \rangle \mapsto \sigma _{(\langle \eta ,z_0 \rangle ,\mu _\nu )}(i)$$
is cofinal in $\nu$. If $n_\nu > 1$, then $F:\gamma \times \omega \rightarrow \alpha _{\tau (n_\nu -1, \mu _\nu )}$ where 
$$\langle \eta ,i \rangle \mapsto h_{\alpha _{\tau (n_\nu -1, \mu _\nu )},K^{n_\nu -1}_{\mu _\nu}}(i, \langle \eta ,z_0 ^\ast \rangle )$$
is cofinal in $\nu$ by the proof of Lemma 12 (vi), where
\smallskip 

$f^{n_\nu -1}_{(\beta ,\langle x (n_\nu -1,\mu _\nu ),z_0 \rangle ,\mu _\nu )}(z_0 ^\ast )=z_0$.
\smallskip \\
But $F$ is definable over $I_{\mu _\nu}$ by lemma 13. On the other hand, in a standard morass,
\smallskip 

$\nu \vartriangleleft \tau$ $\Rightarrow$ $\nu$ regular in $J^D_\tau$.
\smallskip \\
So $\nu$ is regular in $I_{\mu _\nu}$.
$\Box$
\medskip \\
{\bf Remark}
\smallskip \\
Let $L[X]$ satisfy (Amenability), (Coherence) and (Condensation). Then we can do fine structure arguments, especially we have the $\Sigma _n$-Skolem functions $h^n_\nu$ of $I_\nu$. As a result, we get: If $S^X =\{ \beta (\nu )\mid \nu$ singular in $I_\kappa \}$, then $S^X =\{ \nu \mid \nu$ singular in $I_{\nu +\omega} \}$. Because $\{ \nu \mid \nu$ singular in $I_{\nu +\omega}\} \subseteq \{ \beta (\nu )\mid \nu$ singular in $I_\kappa \}$ by definition. For $\{ \beta (\nu )\mid \nu$ singular in $I_\kappa \} \subseteq \{ \nu \mid \nu$ singular in $I_{\nu +\omega}\}$, let $n$ be least such that $\nu$ becomes singular over $I_{\mu _\nu}$. Let $p$ be minimal such that $\nu$  becomes singular over $I_{\mu _\nu}$ in the parameter $p$. Let $p^\ast$ be minimal such that  $h^n_{\mu _\nu}(i,p^\ast )=p$ for some $i \in \omega$. Let $\pi :I_{\bar \mu} \rightarrow I_{\mu _\nu}$ be the uncollapse of $h^n_{\mu _\nu}[\omega \times (J^X_\nu \times \{ p^\ast \} )]$. Let $\pi (\bar p)=p^\ast$. Then $\nu$ becomes singular over $I_{\bar \mu}$ and $h^n_{\bar \mu}[\omega \times (J^X_\nu \times \{ \bar p \} )]=J^X_{\bar \mu}$. By the minimality of $\mu _\nu$, we get $\bar \mu =\mu _\nu$ and that $\mu _\nu \in \{ \nu \mid \nu$ is singular in $I_{\nu +\omega}\}$.

Conversely, if $S^X =\{ \nu \mid \nu$ singular in $I_{\nu +\omega}\}$, then $S^X =\{ \beta (\nu )\mid \nu$ singular in $I_\kappa \}$. We prove $\{ \beta (\nu )\mid \nu$ singular in $I_\kappa \} \subseteq \{ \nu \mid \nu$ singular in $I_{\nu +\omega}\}$ as above. And $\{ \nu \mid \nu$ singular in $I_{\nu +\omega}\} \subseteq \{ \beta (\nu )\mid \nu$ singular in $I_\kappa \}$ holds again by definition.

\section*{References}
{[BJW]} A. Beller, R. Jensen, P. Welch: {\bf Coding the universe}, London Mathematical Society Lecture Notes Series, vol. 47, Cambridge University Press, Cambridge, 1982
\smallskip\\
{[ChKe]} C. C. Chang, H. J. Keisler: {\bf Model Theory}, North-Holland, Amsterdam, 1999
\smallskip\\
{[Dev]} K. Devlin: {\bf Constructibility}, Springer-Verlag, Berlin, 1984
\smallskip\\
{[Don]} D. Donder: {\it Another look at gap-1 morasses}, {\bf Recursion theory}, Proceedings of symposia in pure mathematics, vol. 42, American mathematical society, Providence, RI, 1985, 223 - 236
\smallskip\\
{[DJS]} H.-D. Donder, R. B. Jensen, L. J. Stanley: {\it Condensation-Coherent Global Square Systems}, {\bf Proceedings of Symposia in Pure Mathematics}, vol. 42, 1985
\smallskip\\
{[Fr]} S. D. Friedman: {\bf Fine Structure and Class Forcing}, De Gruyter, Berlin, 2000
\smallskip\\
{[Irr1]} B. Irrgang: {\bf Kondensation und Moraste}, Dissertation, M\"unchen, 2002
\smallskip\\
{[Irr2]} B. Irrgang: {\it Constructing $(\omega _1, \beta)$-morasses, $\beta \geq \omega _1$}
\smallskip\\
{[Jen]} R. Jensen: {\bf Higher-Gap Morasses}, hand-written notes, 1972/73
\smallskip\\
{[JeZe]} R. Jensen, M. Zeman: {\it Smooth categories and global $\Box$}, {\bf Annals of Pure and Applied Logic} 102 (2000), 101 - 138
\smallskip\\
{[Mor]} C. Morgan: {\it Higher gap morasses, Ia: Gap-two morasses and condensation}, {\bf The journal of symbolic logic}, vol. 63, no. 3, 1998, 753 - 787
\smallskip\\
{[SchZe]} R. Schindler, M. Zeman: {\it Fine Structure Theory}, to appear in {\bf Handbook of Symbolic Logic}, edited by M. Foreman, A. Kanamori, M. Magidor
\smallskip\\
{[SchZe1]} E. Schimmerling, M. Zeman: {\it Characterisation of $\Box _\kappa$ in Core Models}, {\bf Journal of Mathematical Logic} 4 (2004), 1 - 72
\smallskip\\
{[SchZe2]} E. Schimmerling, M. Zeman: {\it Cardinal Transfer Properties in Extender Models I}
\smallskip\\
{[She]} S. Shelah: {\bf Cardinal Arithmetic}, Clarendon Press, Oxford, 1994 
\smallskip\\
{[ShSt1]} S. Shelah, L. Stanley: {\it S-forcing I: A ``black box'' theorem for morasses, with applications to super-Souslin trees}, {\bf Israel Journal of Mathematics}, vol. 43, no. 3, 1982, 185 - 236 
\smallskip\\
{[ShSt2]} S. Shelah, L. Stanley: {\it The combinatorics of combinatorial coding by a real}, {\bf The Journal of symbolic logic} 60 (1995), 36 - 57 
\smallskip\\
{[ShSt3]} S. Shelah, L. Stanley: {\it A combinatorial forcing for coding the universe by a real when there are no sharps}, {\bf The Journal of symbolic logic} 60 (1994), 1 - 35
\smallskip\\
{[Sta1]} L. Stanley: {\bf L-like models of set theory: Forcing, combinatorial principles, and morasses}, Dissertation, UC Berkeley, 1977 
\smallskip\\
{[Sta2]} L. Stanley: {\it A short course on gap-one morasses with a review of the fine structure of L}, {\bf Surveys in Set Theory}, London Mathematical Society Lecture Notes Series, vol. 87, Cambridge University Press, Cambridge 1983, 197 - 243
\smallskip\\
{[Vel1]} D. Velleman: {\it Simplified Morasses}, {\bf The Journal of Symbolic Logic}, vol. 49, no. 1, 1984, 257 - 271 
\smallskip\\
{[Vel2]} D. Velleman: {\it Souslin trees constructed from morasses}, {\bf Axiomatic Set Theory (Boulder, Colo., 1983)}, Contemporary Mathematics, vol. 31, Amer. Math. Soc., Providence, RI, 1984, 219 - 241 
\smallskip\\
{[Vel3]} D. Vellman: {\it Simplified morasses with linear limits}, {\bf The Journal of symbolic logic}, vol. 49, no. 4, 1984, 1001 - 1021
\smallskip\\
{[Vel4]} D. Velleman: {\it Simplified gap-2 morasses}, {\bf Annals of Pure and Applied Logic}, vol. 34, 1987, 171 - 208 
\end{document}